\documentclass[reqno,11pt]{amsart}
\DeclareMathOperator*{\esssup}{\mathrm{ess\,sup}}

\usepackage{amsmath,amssymb}
\usepackage{dsfont}         

\newtheorem{theorem}{Theorem}[section]
\newtheorem{corollary}[theorem]{Corollary}
\newtheorem{lemma}[theorem]{Lemma}
\newtheorem{proposition}[theorem]{Proposition}
\theoremstyle{definition}
\newtheorem{definition}[theorem]{Definition}
\theoremstyle{remark}
\newtheorem{remark}[theorem]{Remark}
\numberwithin{equation}{section}

\title{A Widom-Rowlinson Jump Dynamics in the Continuum}

\author{Joanna Bara{\'n}ska}

\address{Instytut Matematyki, Uniwersytet Marii Curie-Sk{\l}odowskiej, 20-031 Lublin, Poland}
\email{asia13p@wp.pl}

\author{ Yuri  Kozitsky}

\address{Instytut Matematyki, Uniwersytet Marii Curie-Sk{\l}odowskiej, 20-031 Lublin, Poland}
\email{jkozi@hektor.umcs.lublin.pl}

\begin{document}



\begin{abstract}
We study the dynamics of an infinite system of point particles of
two types. They perform random jumps in $\mathds{R}^d$ in the course
of which particles of different types repel each other whereas those
of the same type do not interact. The states of the system are
probability measures on the corresponding configuration space, the
global (in time) evolution of which is constructed by means of
correlation functions. It is proved that for each initial
sub-Poissonian state $\mu_0$, the states evolve $\mu_0 \mapsto
\mu_t$ in such a way that $\mu_t$ is sub-Poissonian for all $t>0$.
The mesoscopic (approximate) description of the evolution of states
is also given. The stability of translation invariant stationary
states is studied. In particular, we show that some of such states
can be unstable with respect to space-dependent perturbations.

\end{abstract}
\maketitle


\section{Introduction}
\label{S1}

\subsection{Posing}
In this paper, we study the dynamics of an infinite system of point
particles of two types placed in $\mathds{R}^d$. The particles
perform random jumps in the course of which particles of different
types repel each other whereas those of the same type do not
interact. We do not require that the repulsion is of hard-core type.
This model can be viewed as a dynamical version of the
Widom-Rowlinson model \cite{WR} of equilibrium  statistical
mechanics  -- one of the few models of phase transitions in
continuum particle systems, see the corresponding discussion in
\cite{FKKO} where a similar birth-and-death model was introduced and
studied.

The phase space of our model is defined as follows. Let $\Gamma$
denote the set of all  $\gamma \subset \mathds{R}^d$ that are {\it
locally finite}, i.e., such that $\gamma\cap\Lambda$ is a finite set
whenever $\Lambda \subset \mathds{R}^d$ is compact. Thus, $\Gamma$
is a configuration space as defined in \cite{Albev,BKKK,FKO,Tobi}.
In order to take into account the particle's type we use the
Cartesian product $\Gamma^2 = \Gamma \times \Gamma$, see
\cite{F,FKKO,FKO1}, the elements of which are denoted by $\gamma =
(\gamma_0, \gamma_1)$. In a standard way, $\Gamma^2$ is equipped
with a $\sigma$-field of measurable subsets which allows one to deal
with probability measures considered as states of the system. Among
them one may distinguish Poissonian states in which the particles
are independently distributed over $\mathds{R}^d$. In
\emph{sub-Poissonian} states, the dependence between the particle's
positions  is not too strong. As was shown in \cite{KK}, for
infinite particle systems with birth-and-death dynamics the states
remain sub-Poissonian globally in time if the birth of the particles
is in a sense controlled by their state-dependent death. In
\cite{FKKO}, the evolution of sub-Poissonian correlation function of
a birth-and-death Widom-Rowlinson model was shown to hold on a
bounded time interval. For conservative dynamics in which the
particles just change their positions, the interaction may in
general change the sub-Poissonian character of the state in finite
time (even cause an explosion), e.g., due to an infinite number of
simultaneous correlated jumps.   Thus, the conceptual outcome of the
present study is that this is not the case for the considered model.
Our another aim in this paper is to study the dynamics of the
considered model in the mesoscopic limit, which yields its though an
approximate (mean-field like) but more detailed picture.

\subsection{Presenting the results}

The evolution of systems like the one we consider is described by
the Kolmogorov equation
\begin{equation}
  \label{1}
 \frac{d}{dt} F_t = L F_t, \qquad F_t|_{t=0} = F_0,
\end{equation}
where $F_t:\Gamma^2\to\mathds{R}$ is an {\it observable} and the
operator $L$ specifies the model. In our case it has the following
form
\begin{eqnarray}\label{LF}
(LF)(\gamma_{0},\gamma_{1}) &= &\sum_{x\in
\gamma_{0}}\int_{\mathds{R}^d} c_{0}(x,y,\gamma_{1})
[F(\gamma_{0}\backslash x \cup
y,\gamma_{1})-F(\gamma_{0},\gamma_{1})]d y \qquad \\[.2cm] &+ & \sum_{x\in
\gamma_{1}}\int_{\mathds{R}^d} c_{1}(x,y,\gamma_{0})
[F(\gamma_{0},\gamma_{1}\backslash x \cup
y)-F(\gamma_{0},\gamma_{1})]d y. \nonumber
\end{eqnarray}
The evolution of states is supposed to be obtained by solving the
Fokker-Planck equation
\begin{equation}
  \label{1a}
\frac{d}{dt} \mu_t = L^* \mu_t, \qquad \mu_t|_{t=0} = \mu_0,
\end{equation}
related to that in (\ref{1}) by the duality
\begin{equation}
  \label{1b}
  \int_{\Gamma^2}F_t (\gamma) \mu_0 ( d \gamma) = \int_{\Gamma^2}F_0 (\gamma) \mu_t ( d
  \gamma).
\end{equation}
As is usual for models of this kind, the direct meaning of (\ref{1})
or (\ref{1a}) can only be given for states of finite systems, cf.
\cite{K}. In this case, the Banach space where the Cauchy problem in
(\ref{1a}) is defined can be the space of signed measures with
finite variation. For infinite systems, the evolution of states is
constructed by means of correlation functions, see
\cite{BKKK,FKKK,FKKO,FKO,FKO1,KK} and the references quoted therein.

In this paper, in describing the evolution of states, see Theorem
\ref{1tm} below,  we mostly follow the scheme elaborated in
\cite{KK}. It consists in: (a) constructing the evolution of
correlation functions $k_0 \mapsto k_t$, $t< T < +\infty$, based on
the Cauchy problem in (\ref{16}); (b) proving that each $k_t$ is the
correlation function of a unique sub-Poissonian state $\mu_t$; (c)
constructing the continuation of thus obtained evolution
$k_{\mu_0}=k_0 \mapsto k_t = k_{\mu_t}$ to all $t>0$. Step (a) is
performed by means of Ovcyannikov-like arguments similar to those
used, e.g., in \cite{BKKK,FKKK,FKKO}. Step (b) is based on the use
of the Denjoy-Carleman theorem \cite{DC}. In realizing step (c), we
crucially use the result of (b). Our description of the mesoscopic
limit is based on the a scaling procedure described in Section
\ref{S4}. It is equivalent to the Lebowitz-Penrose scaling used in
\cite{FKKO}, and also to the Vlasov scaling used in
\cite{BKKK,FKKK}. In this procedure, passing to the mesoscopic level
amounts to considering the system at different spatial scales
parameterized by $\varepsilon \in (0; 1]$ in such a way that
$\varepsilon = 1$ corresponds to the micro-level, whereas the limit
$\varepsilon \to 0$ yields the meso-level description in which the
corpuscular structure disappears and the system turns into a
(two-component) medium characterized by a density function. The
evolution of the latter is supposed to be found from the kinetic
equation (\ref{29}). In Theorem \ref{Ktm}, we show that the kinetic
equation has a unique global (in time) solution in the corresponding
Banach space. In Theorem \ref{4tm}, we demonstrate that the micro-
and mesoscopic descriptions are indeed connected by the scaling
procedure in the sense of Definition \ref{Ja1df}. In Theorems
\ref{K1tm} and \ref{K2tm}, we describe the stability of translation
invariant stationary solutions of the kinetic equation. In
particular, we show that some of such solutions can be unstable with
respect to space-dependent perturbations.

The rest of the paper has the following structure. In Section
\ref{Sq2}, we give necessary information on the analysis in two
component configuration spaces and on the description of
sub-Poissonian states on such spaces with the help of Bogoliubov
functionals and correlation functions. We also describe in detail
the model which we consider. In Section \ref{S3}, we formulate the
results mentioned above and prove Theorems \ref{K1tm} and
\ref{K2tm}. We also provide some comments; in particular, we relate
our results with those of \cite{FKKO} describing a birth-and-death
version of the Widom-Rowlinson dynamics in the continuum. Section
\ref{S4} is dedicated to developing our main technical tool --
Proposition \ref{3tm}. By means of it we realize step (a) in proving
Theorem \ref{1tm}, see above. Steps (b) and (c) are  based on Lemmas
\ref{Id1lm}, \ref{A1lm}, \ref{J1lm} and \ref{J2lm} proved in Section
\ref{SPr}. Section \ref{S6} is dedicated to the proof of Theorems
\ref{Ktm} and \ref{K1tm}.

\section{Preliminaries and the Model}
\label{Sq2}
\subsection{Two-component configuration spaces}

Here we briefly present necessary information on the subject. A more
detailed description can be found in, e.g.,  \cite{F,FKKO,FKO1}.

Let $\mathcal{B}(\mathds{R}^d)$ and $\mathcal{B}_{\rm
b}(\mathds{R}^d)$ denote the sets of all Borel and all bounded Borel
subsets of $\mathds{R}^d$, respectively. The configuration space
$\Gamma$ mentioned above is equipped with the vague topology and
thus with the corresponding Borel $\sigma$-field
$\mathcal{B}(\Gamma)$. Then $\Gamma^2:= \Gamma \times \Gamma$ is
equipped with the product $\sigma$-field which we denote by
$\mathcal{B}(\Gamma^2)$. The elements of $\Gamma^2$ are $\gamma =
(\gamma_0, \gamma_1)$, i.e., the one-component configurations are
always written with the subscript $i=0,1$. Likewise, for
$\Lambda_i\in \mathcal{B}(\mathds{R}^d)$, $i=0,1$, we denote
$\Lambda = \Lambda_0 \times \Lambda_1$ and set
\[
\Gamma^2_\Lambda = \{\gamma=(\gamma_0, \gamma_1)\in \Gamma^2:
\gamma_i \subset \Lambda_i, \ i=0,1\}.
\]
Clearly $\Gamma^2_\Lambda \in \mathcal{B}(\Gamma^2)$ and hence
\[
 \mathcal{B}(\Gamma^2_\Lambda):=\{ A \cap \Gamma^2_\Lambda : A \in \mathcal{B}(\Gamma^2)\}
\]
is a sub-field of $\mathcal{B}(\Gamma^2)$. Let
$p_{\Lambda}:\Gamma^2\to \Gamma^2_\Lambda$ be the projection
\begin{equation*}
 p_\Lambda (\gamma) = \gamma_\Lambda := (\gamma_0 \cap \Lambda_0, \gamma_1 \cap \Lambda_1).
\end{equation*}
It is clearly measurable, and thus the sets
\begin{equation*}
 p^{-1}_\Lambda(A_\Lambda) :=\{ \gamma\in \Gamma^2: p_\Lambda (\gamma) \in A_\Lambda \},
 \quad A_\Lambda \in \mathcal{B}(\Gamma^2_\Lambda),
\end{equation*}
belong to $\mathcal{B}(\Gamma^2)$ for each Borel $\Lambda_i$,
$i=0,1$.

Let $\mathcal{P}(\Gamma^2)$ denote the set of all  probability
measures on $(\Gamma^2, \mathcal{B}(\Gamma^2))$. For a given $\mu\in
\mathcal{P}(\Gamma^2)$, its projection on $(\Gamma^2_\Lambda,
\mathcal{B} (\Gamma^2_\Lambda))$ is
\begin{equation}
 \label{5}
\mu^\Lambda (A_\Lambda) := \mu\left(p^{-1}_\Lambda (A_\Lambda)
\right), \qquad A_\Lambda \in \mathcal{B}(\Gamma^2_\Lambda).
\end{equation}
Let $\pi$ be the standard homogeneous Poisson measure on
$(\Gamma,\mathcal{B}(\Gamma))$ with density (intensity) $\varkappa
=1$. Then the product measure $\pi^2:=\pi\otimes \pi$ is a
probability measure on $(\Gamma^2, \mathcal{B}(\Gamma^2))$. By
$\mathcal{P}_\pi (\Gamma^2)$ we denote the set of all $\mu\in
\mathcal{P}(\Gamma^2)$, for each of which the projections
$\mu^\Lambda$, with all possible $\Lambda = \Lambda_0 \times
\Lambda_1$, $\Lambda_i \in \mathcal{B}_{\rm b} (\mathds{R}^d)$,
$i=0,1$, are absolutely continuous with respect to the corresponding
projections of $\pi^2$. It is known, see \cite[Proposition 3.1]{F},
that for each $\mu \in \mathcal{P}_\pi (\Gamma^2)$ the following
holds
\begin{equation*}
\mu\left(\{ \gamma=(\gamma_0, \gamma_1) \in \Gamma^2: \gamma_0 \cap
\gamma_1 = \emptyset\}\right) =1.
\end{equation*}
Since we are going to deal with elements of $\mathcal{P}_\pi
(\Gamma^2)$ only, from now on we assume that the configurations
$\gamma_0$ and $\gamma_1$ are subsets of different copies of
$\mathds{R}^d$.

 Let $\Gamma_0^2$ be the set of all finite $\gamma \in
\Gamma^2$. It is an element of $\mathcal{B}(\Gamma^2)$ as each of
$\gamma \in \Gamma_0^2$ belongs to a certain $\Gamma^2_\Lambda$,
$\Lambda = \Lambda_0\times \Lambda_1$, $\Lambda_{i} \in
\mathcal{B}_{\rm b}(\mathds{R}^d)$, $i=0,1$. Note that
$\Gamma^2_\Lambda \subset \Gamma^2_0$ for each such $\Lambda$. Set
$\mathds{N}_0 = \mathds{N}\cup \{0\}$, and then also $\mathds{N}_0^2
= \{n=(n_0,n_1): n_i \in \mathds{N}_0, i=0,1\}$. It can be proved
that a function $G:\Gamma^2_0 \to \mathds{R}$ is
$\mathcal{B}(\Gamma^2)/\mathcal{B}(\mathds{R} )$-measurable if and
only if for each $n\in \mathds{N}_0^2$, there exists a Borel
function $G^{(n)}: (\mathds{R}^{d})^{n_0} \times
(\mathds{R}^{d})^{n_1} \to \mathds{R}$, symmetric with respect to
the permutations of the components of each of $\eta_i$, $i=0,1$,
such that
\begin{equation*}
 G(\eta) = G(\eta_0, \eta_1) = G^{(n)} ( x_1, \dots , x_{n_0};  y_1, \dots , y_{n_1}),
\end{equation*}
for $\eta_0 = \{ x_1, \dots , x_{n_0}\}$ and $\eta_1 =  \{y_1, \dots
, y_{n_1}\}$.
\begin{definition}
  \label{Gdef}
By $B_{\rm bs}(\Gamma_0^2)$ we denote the set of all measurable
functions $G:\Gamma^2_0 \to \mathds{R}$ that have the following two
properties: (a) there exists $\Lambda = \Lambda_0 \times \Lambda_1$
with $\Lambda_i \in \mathcal{B}_{\rm b} (\mathds{R}^d)$, $i=0,1$,
such that $G(\eta) = 0$ whenever $\eta_i\cap \Lambda_i^c \neq
\emptyset$ for either of $i=0,1$; (b) there exists $N\in
\mathds{N}_0$ such that $G(\eta)=0$ whenever $\max_{i=0,1}|\eta_i|
>N$. Here $\Lambda_i^c := \mathds{R}^d
\setminus \Lambda_i$ and $|\cdot |$ stands for cardinality. By
$\Lambda(G)$ and $N(G)$ we denote the smallest $\Lambda$ and $N$
with the properties just described.
\end{definition}
The Lebesgue-Poisson measure $\lambda$ on $(\Gamma^2_0,
\mathcal{B}(\Gamma^2_0))$ is then defined by the following formula
\begin{eqnarray}
\label{8}
& & \int_{\Gamma_0^2} G(\eta ) \lambda ( d \eta) \\[.2cm]
& & \quad = \sum_{n_0=0}^\infty \sum_{n_1=0}^\infty \frac{1}{n_0! n_1!} \int_{(\mathds{R}^d)^{n_0}}
\int_{(\mathds{R}^d)^{n_1}} G^{(n)} ( x_1, \dots , x_{n_0};  y_1, \dots , y_{n_1}) \nonumber \\[.2cm]
& & \quad \times d x_1 \cdots dx_{n_0} d y_1 \cdots dy_{n_1}, \nonumber
\end{eqnarray}
which has to hold for all $G\in B_{\rm bs}(\Gamma_0^2)$ with the
usual convention regarding the cases $n_i=0$. The same can also be
written as
\begin{equation}
 \label{9}
\int_{\Gamma_0^2} G(\eta ) \lambda ( d \eta)  = \int_{\Gamma_0}
\int_{\Gamma_0} G(\eta_0,\eta_1 ) (\lambda_0 \otimes \lambda_1)( d
\eta_0, d \eta_1)
\end{equation}
where both $\lambda_i$ are the copies of the standard
Lebesgue-Poisson measure on the single-component set $\Gamma_0$. In
the sequel, both Lebesgue-Poisson measures on $\Gamma_0^2$ and on
$\Gamma_0$ will be denoted by $\lambda$ if no ambiguity may arise.

For $\gamma\in \Gamma^2$, by writing $\eta \Subset \gamma$ we mean
that $\eta_i\Subset \gamma_i$, $i=0,1$, i.e., both $\eta_i$ are
nonempty and finite. For $G\in B_{\rm bs}(\Gamma^2)$, we set
\begin{equation}
  \label{9a}
(KG)(\gamma) := \sum_{\eta \Subset \gamma} G(\eta) = \sum_{\eta_0
\Subset \gamma_0} \sum_{\eta_1 \Subset \gamma_1} G(\eta_0,\eta_1).
\end{equation}
Note that the sums in (\ref{9a}) are finite and $KG$ is a cylinder
function on $\Gamma^2$. The latter means that it is
$\mathcal{B}(\Gamma^2_{\Lambda(G)})$-measurable, see Definition
\ref{Gdef}. Moreover,
\begin{equation}
  \label{9b}
|(KG)(\gamma)| \leq \left( 1 + |\gamma_0\cap\Lambda_0
(G)|\right)^{N_0(G)} \left( 1 + |\gamma_1\cap\Lambda_1
(G)|\right)^{N_1(G)}.
\end{equation}

\subsection{Correlation functions}

In the approach we follow, see \cite{BKKK,FKKK,KK}, the evolution of
states is constructed in the next way. Let $\varTheta$ denote the
set of all compactly supported continuous maps $\theta = (\theta_0,
\theta_1): \mathds{R}^d \times \mathds{R}^d \to (-1,0]^2$. For each
$\theta \in \varTheta$, the map
$$\Gamma^2 \ni \gamma \mapsto \prod_{x\in \gamma_0} (1+ \theta_0
(x)) \prod_{y\in \gamma_1} (1+ \theta_1 (y))$$ is measurable and
bounded. Hence, for a state $\mu$, one may define
\begin{equation}
  \label{1d}
  B_\mu (\theta) = \int_{\Gamma^2} \prod_{x\in \gamma_0} (1+ \theta_0
(x)) \prod_{y\in \gamma_1} (1+ \theta_1 (y)) \mu ( d \gamma),
\end{equation}
-- the so called  Bogoliubov functional for $\mu$, considered as a
map $\varTheta \to \mathds{R}$. Let $\mathcal{P}_{\rm
exp}(\Gamma^2)$ stand for the set of $\mu \in \mathcal{P}_\pi
(\Gamma^2)$ for which $B_\mu$ can be extended to an exponential type
entire function of $\theta \in L^1 (\mathds{R}^d\times
\mathds{R}^d\to \mathds{R}^2) $. This exactly means that $B_\mu$ can
be written in the form, cf. (\ref{9}),
\begin{equation}
  \label{1e}
 B_\mu (\theta) = \int_{\Gamma_0^2} k_\mu (\eta) E(\eta;\theta)
 \lambda (d \eta),
\end{equation}
cf. (\ref{8}), with $k_\mu:\Gamma_0^2 \to [0,+\infty)$ such that
$k_\mu^{(n)}\in L^\infty ((\mathds{R}^d)^{n_0}\times
(\mathds{R}^d)^{n_1}\to \mathds{R})$ and
\begin{equation}
  \label{1f}
 E(\eta;\theta) = e (\eta_0; \theta_0)e (\eta_1; \theta_1)  := \prod_{x\in \eta_0} \theta_0 (x) \prod_{y\in \eta_1} \theta_1
 (y).
\end{equation}
This, in particular, means that $k_\mu$ is essentially bounded with
respect to the Lebesgue-Poisson measure $\lambda$ defined in
(\ref{8}). For the (heterogeneous) Poisson measure $\pi_\varrho$,
the Bogoliubov functional is
\begin{equation}
  \label{1fa}
B_{\pi_\varrho} (\theta) = \exp\left(\sum_{i=0,1}\int_{\mathds{R}^d}
\theta_i (x) \varrho_i(x) dx \right),
\end{equation}
where $\varrho = (\varrho_0, \varrho_1)$ is the (two-component)
density function. Then by (\ref{8}) and (\ref{1e}) we have
\begin{equation}
  \label{1fb}
 k_{\pi_\varrho} (\eta)=E(\eta; \varrho) = e(\eta_0, \varrho_0)e(\eta_1,
 \varrho_1).
\end{equation}
If one rewrites (\ref{1d}) in the form
\[
B_\mu (\theta) = \int_{\Gamma^2} F_\theta (\gamma) \mu(d \gamma),
\]
then the action of $L$ on $F$ as in (\ref{LF}) can be transformed to
the action of $L^\Delta$ on $k_\mu$ from the following relation
\begin{equation}
  \label{1g}
\int_{\Gamma^2}(L F_\theta) (\gamma) \mu(d \gamma) =
\int_{\Gamma_0^2} (L^\Delta k_\mu) (\eta) E(\eta;\theta)\lambda (d
\eta)
\end{equation}
The main advantage of this is that $k_\mu$ is  a function of {\em
finite} configurations.

For $\mu \in \mathcal{P}_{\rm exp}(\Gamma^2)$ and $\Lambda
=(\Lambda_0, \Lambda_1)$, $\Lambda_i \in \mathcal{B}_{\rm
b}(\mathds{R}^d)$, let $\mu^\Lambda$ be as in (\ref{5}). Then
$\mu^\Lambda$ is absolutely continuous with respect to the
corresponding restriction  $\lambda^\Lambda$ of the measure defined
in (\ref{8}), and hence we may write
\begin{equation}
\label{9c} \mu^\Lambda (d \eta ) = R^\Lambda_\mu (\eta)
\lambda^\Lambda ( d \eta), \qquad \eta \in \Gamma^2_\Lambda.
\end{equation}
Then the correlation function $k_\mu$ and the Radon-Nikodym
derivative $R_\mu^\Lambda$ are related to each other by, cf.
(\ref{9}),
\begin{eqnarray}
  \label{9d}
k_\mu(\eta) & = & \int_{\Gamma^2_\Lambda} R^\Lambda_\mu (\eta \cup
\xi) \lambda^\Lambda ( d\xi)\\[.2cm]& = & \int_{\Gamma_{\Lambda_0}}
\int_{\Gamma_{\Lambda_1}} R^\Lambda_\mu (\eta_0 \cup\xi_0, \eta_1
\cup \xi_1) (\lambda_0^{\Lambda_0}  \otimes \lambda_1^{\Lambda_1})
(d\xi_0, d \xi_1), \qquad \eta \in \Gamma^2_\Lambda. \nonumber
\end{eqnarray}
Note that (\ref{9d}) relates $R^\Lambda_\mu$ with the restriction of
$k_\mu$ to $\Gamma_\Lambda^2$. The fact that these are the
restrictions of one and the same function
$k_\mu:\Gamma_0^2\to\mathds{R}$ corresponds to the Kolmogorov
consistency of the family $\{\mu^\Lambda\}_{\Lambda}$.

By (\ref{9a}), (\ref{5}), and (\ref{9c}) we get
\begin{equation*}
\int_{\Gamma^2} (KG)(\gamma) \mu(d\gamma) = \langle \! \langle G,
k_\mu \rangle \!\rangle,
\end{equation*}
holding for each $G\in B_{\rm bs}(\Gamma_0^2)$ and $\mu \in
\mathcal{P}_{\rm exp}(\Gamma^2)$. Here
\begin{equation}
  \label{9f}
\langle \! \langle G, k \rangle \!\rangle := \int_{\Gamma_0^2}
G(\eta) k(\eta) \lambda (d \eta),
\end{equation}
for suitable $G$ and $k$. Define
\begin{equation}
  \label{9g}
B^\star_{\rm bs} (\Gamma_0^2) =\{ G\in B_{\rm bs}(\Gamma_0^2):
(KG)(\gamma) \geq 0 \ {\rm for} \ {\rm all} \ \gamma\in \Gamma^2\}.
\end{equation}
By \cite[Theorems 6.1 and 6.2 and Remark 6.3]{Tobi} one can prove
that the following holds.
\begin{proposition}
  \label{Gpn}
Let  a measurable function $k : \Gamma_0^2 \to \mathds{R}$  have the
following properties:
\begin{eqnarray}
  \label{9h}
& (a) & \ \langle \! \langle G, k \rangle \!\rangle \geq 0, \qquad
{\rm for} \ {\rm all} \ G\in B^\star_{\rm bs} (\Gamma_0^2);\\[.2cm]
& (b) & \ k(\emptyset, \emptyset) = 1; \qquad (c) \ \ k(\eta) \leq
 C^{|\eta_0|+|\eta_1|} ,
\nonumber
\end{eqnarray}
with (c) holding for some $C >0$ and $\lambda$-almost all $\eta\in
\Gamma_0^2$. Then there exists a unique ½$\mu \in \mathcal{P}_{\rm
exp}(\Gamma^2)$ for which $k$ is the correlation function.
\end{proposition}

\subsection{The model}

The model we consider is specified by the operator $L$ given in
(\ref{LF}) where the coefficients are supposed to be of the
following form
\begin{eqnarray}\label{c}
c_{0}(x,y,\gamma_{1}) & = & a_{0}(x-y)\exp\left(-\sum_{z \in \gamma_{1}}\phi_{0}(y-z)\right),\\[.2cm]
c_{1}(x,y,\gamma_{0}) & = & a_{1}(x-y)\exp\left(-\sum_{z \in \gamma_{0}}\phi_{1}(y-z)\right), \nonumber
\end{eqnarray}
with jump kernels $a_{i}: \mathds{R}^d \rightarrow [0,+\infty)$ such
that $a_{i}(x)=a_{i}(-x)$ and
\begin{equation}
 \label{10}
\int_{\mathds{R}^d}a_{i}(x)d x =: \alpha_i < \infty, \qquad i=0,1.
\end{equation}
The repulsion potentials in (\ref{c}) $\phi_{i}: \mathds{R}^d
\rightarrow [0,+\infty)$ are supposed to be symmetric, $\phi_i(x) =
\phi_i(-x)$, and such that
\begin{eqnarray}
  \label{33}
\int_{\mathds{R}^d} \phi_i(x) d x =:\langle \phi_i \rangle < \infty,
\qquad \esssup_{x\in \mathds{R}^d} \phi_i (x) =:\bar{\phi}_i <
\infty.
\end{eqnarray}
Then
\begin{equation}
 \label{11}
 \int_{\mathds{R}^d}\bigg{(}1-\exp(-\phi_{i}(x))\bigg{)}d x \leq \langle \phi_{i} \rangle, \qquad i=0,1.
\end{equation}
By (\ref{LF}) and (\ref{1g}) one obtains the action of $L^\Delta$ in
the following form. For $x\in \mathds{R}^d$, we set
\begin{equation}
  \label{12}
 \tau_x^i (y) = \exp(-\phi_i (x-y)), \quad  t_x^i (y)= \tau_x^i (y) - 1, \quad
 y\in \mathds{R}^d, \ \ i=0,1.
\end{equation}
Next, for a function $k(\eta) = k(\eta_0, \eta_1)$, cf. (\ref{9}),
we introduce the maps
\begin{eqnarray}
  \label{13}
(Q_y^0 k) (\eta_0, \eta_1) & = & \int_{\Gamma_0} k(\eta_0, \eta_1
\cup \xi)e(t^0_y ;\xi) \lambda(d\xi), \\[.2cm]
(Q_y^1 k) (\eta_0, \eta_1) & = & \int_{\Gamma_0} k(\eta_0\cup \xi,
\eta_1 )e(t^1_y ;\xi) \lambda(d\xi), \nonumber
\end{eqnarray}
where $e$ is as in (\ref{1f}). Then
\begin{eqnarray}
  \label{15}
(L^\Delta k) (\eta_0, \eta_1) & = & \sum_{y\in \eta_0}
\int_{\mathds{R}^d} a_0 (x-y) e(\tau^0_y;\eta_1) (Q_y^0 k)
(\eta_0\setminus y \cup x,
\eta_1) d x \nonumber \\[.2cm]
& - & \sum_{x\in \eta_0} \int_{\mathds{R}^d} a_0 (x-y)
e(\tau^0_y;\eta_1) (Q_y^0 k) (\eta_0,
\eta_1) d y  \\[.2cm]
& + & \sum_{y\in \eta_1} \int_{\mathds{R}^d} a_1 (x-y)
e(\tau^1_y;\eta_0) (Q_y^1 k) (\eta_0,
\eta_1\setminus y \cup x) d x \nonumber \\[.2cm] & - & \sum_{x\in \eta_1} \int_{\mathds{R}^d} a_1 (x-y)
e(\tau^1_y;\eta_0) (Q_y^1 k) (\eta_0, \eta_1) d y. \nonumber
\end{eqnarray}

\section{The results}

 \label{S3}

\subsection{The microscopic level}

As mentioned above, instead of directly studying the evolution of
states by solving the problem in (\ref{1a}), we pass from $\mu_0$ to
the corresponding correlation function $k_{\mu_0}$ and then consider
the problem
\begin{equation}
  \label{16}
\frac{d}{dt} k_t = L^\Delta k_t, \qquad k_t|_{t=0} = k_{\mu_0},
\end{equation}
where $L^\Delta$ is given in (\ref{15}). For this problem, we prove
the existence of a unique global solution $k_t$ which is the
correlation function of a unique state $\mu_t \in \mathcal{P}_{\rm
exp}(\Gamma^2)$.

We begin by defining the problem (\ref{16}) in the corresponding
spaces of functions $k:\Gamma_0^2 \to \mathds{R}$. From the very
representation (\ref{1e}), see also (\ref{8}), it follows that $\mu
\in \mathcal{P}_{\rm exp}(\Gamma^2)$ implies
\begin{equation*}
 |k_\mu (\eta)| \leq C \exp\bigg{(} \vartheta  \left( |\eta_0| +
|\eta_1|\right)\bigg{)},
\end{equation*}
holding for $\lambda$-almost all $\eta\in \Gamma_0^2$, some $C>0$,
and $\vartheta\in \mathds{R}$. Keeping this in mind we set
\begin{equation}
  \label{17a}
 \|k\|_\vartheta = \esssup_{\eta \in \Gamma^2_0}\left\{ |k_\mu (\eta)| \exp\bigg{(} - \vartheta\left(
  |\eta_0| +|\eta_1| \right)\bigg{)} \right\}.
\end{equation}
Then
\begin{equation*}
\mathcal{K}_\vartheta := \{ k:\Gamma^2_0\to \mathds{R}:
\|k\|_\vartheta <\infty\}
\end{equation*}
is a Banach space with norm (\ref{17a}) and the usual linear
operations. In fact, we are going to use the ascending scale of such
spaces $\mathcal{K}_\vartheta$, $\vartheta \in \mathds{R}$, with the
property
\begin{equation}
  \label{19}
\mathcal{K}_\vartheta \hookrightarrow \mathcal{K}_{\vartheta'},
\qquad \vartheta < \vartheta',
\end{equation}
where $\hookrightarrow$  denotes continuous embedding. Set, cf.
(\ref{9f}) and (\ref{9g}),
\begin{equation}
  \label{19a}
\mathcal{K}^\star_\vartheta =\{k\in \mathcal{K}_\vartheta: \langle
\! \langle G,k \rangle \! \rangle \geq 0 \ {\rm for} \ {\rm all} \
G\in B^\star_{\rm bs} (\Gamma_0^2)\}.
\end{equation}
It is a subset of the cone
\begin{equation}
  \label{19b}
\mathcal{K}^+_\vartheta =\{k\in \mathcal{K}_\vartheta: k(\eta) \geq
0 \ \ {\rm for} \  \lambda-{\rm almost} \ {\rm all} \ \eta \in
\Gamma_0^2\}.
\end{equation}
By Proposition \ref{Gpn} it follows that each $k\in
\mathcal{K}^\star_\vartheta$ such that $k(\emptyset,\emptyset) = 1$
is the correlation function of a unique $\mu\in \mathcal{P}_{\rm
exp}(\Gamma^2)$. Then we define
\begin{equation}
  \label{19c}
\mathcal{K} = \bigcup_{\vartheta \in \mathds{R}}
\mathcal{K}_\vartheta, \qquad \mathcal{K}^\star = \bigcup_{\vartheta
\in \mathds{R}} \mathcal{K}_\vartheta^\star.
\end{equation}
As a sum of Banach spaces, the linear space $\mathcal{K}$ is
equipped with the corresponding inductive topology which turns it
into a locally convex space.

For a given $\vartheta\in \mathds{R}$, by (\ref{12}) -- (\ref{15})
we define $L^\Delta_\vartheta$ as a linear operator in
$\mathcal{K}_\vartheta$ with domain
\begin{equation}
  \label{20}
\mathcal{D} (L^\Delta_\vartheta) = \{ k\in \mathcal{K}_\vartheta:
L^\Delta k \in \mathcal{K}_\vartheta\}.
\end{equation}
\begin{lemma}
  \label{1lm}
For each $\vartheta'' < \vartheta$, cf. (\ref{19}), it follows that
$\mathcal{K}_{\vartheta''} \subset \mathcal{D}
(L^\Delta_\vartheta)$.
\end{lemma}
\begin{proof}
For $\vartheta'' < \vartheta$, by (\ref{11}), (\ref{12}),
(\ref{13}), and (\ref{17a}) we have
\begin{eqnarray}
  \label{21}
\left\vert (Q^0_y k)(\eta_0, \eta_1)\right\vert   & \leq &
\|k\|_{\vartheta''} \exp\left(  \vartheta'' |\eta_0| + \vartheta''
|\eta_1|\right) \qquad  \\[.2cm] & \times & \int_{\Gamma_0} \exp\left(  \vartheta'' |\xi|\right) \prod_{z\in \xi}
\bigg{(} 1 - \exp\left( - \phi_0 (z-y)\right)\bigg{)} \lambda (
d\xi)\nonumber \\[.2cm]
& \leq & \|k\|_{\vartheta''} \exp\left(  \vartheta'' |\eta_0| +
\vartheta'' |\eta_1|\right) \exp\left(\langle \phi_{0} \rangle
e^{\vartheta''} \right) . \nonumber
\end{eqnarray}
Likewise
\begin{equation}
  \label{22}
\left\vert (Q^1_y k)(\eta_0, \eta_1)\right\vert \leq
\|k\|_{\vartheta''} \exp\left( \vartheta'' |\eta_0| + \vartheta''
|\eta_1|\right) \exp\left(\langle \phi_{1} \rangle e^{\vartheta''}
\right) .
\end{equation}
Now we apply the latter two estimates together with (\ref{10}) in
(\ref{15}) and obtain
\begin{eqnarray}
  \label{23}
& & \left\vert (L^\Delta k)(\eta_0, \eta_1)\right\vert \leq 2
\|k\|_{\vartheta''} \exp\left(  \vartheta'' |\eta_0| + \vartheta''
|\eta_1|\right) \\[.2cm]
& & \quad \times \bigg{(} \alpha_0 |\eta_0| \exp\left(\langle
\phi_{0} \rangle e^{\vartheta''} \right) + \alpha_1 |\eta_1|
\exp\left(\langle \phi_{1} \rangle e^{\vartheta''} \right) \bigg{)}.
\nonumber
\end{eqnarray}
By means of the inequality $x\exp(-\sigma x) \leq 1/ e \sigma$, $x,
\sigma
>0$, we get from (\ref{17a}) and (\ref{23}) the following estimate
\begin{gather}
  \label{24}
 \|L^\Delta k\|_{\vartheta}  \leq  \frac{4\|k\|_{\vartheta''}}{e(\vartheta - \vartheta'')}
 \max_{i=0,1}\alpha_i \exp\left(
\langle \phi_{i} \rangle e^{\vartheta''}\right),
\end{gather}
which yields the proof.
\end{proof}
\begin{corollary}
  \label{Gco}
For each $\vartheta,\vartheta''\in \mathds{R}$ such that
$\vartheta'' < \vartheta$, the expression in (\ref{15}) defines a
bounded linear operator $L^\Delta_{\vartheta\vartheta''}:
\mathcal{K}_{\vartheta''}\to \mathcal{K}_{\vartheta}$ the norm of
which can be estimated by means of (\ref{24}).
\end{corollary}
In what follows, we consider two types of operators defined by the
expression in (\ref{15}): (a) unbounded operators
$(L^\Delta_\vartheta, \mathcal{D}(L^\Delta_\vartheta))$,
$\vartheta\in \mathds{R}$, with domains as in (\ref{20}) and Lemma
\ref{1lm}; (b) bounded operators $L^\Delta_{ \vartheta \vartheta''}$
described in Corollary \ref{Gco}. These operators are related to
each other in the following way:
\begin{equation}
  \label{24a}
\forall \vartheta'' < \vartheta \ \  \forall k \in
\mathcal{K}_{\vartheta''} \qquad L^\Delta_{\vartheta
 \vartheta''} k = L^\Delta_{\vartheta} k.
\end{equation}
 By means of the bounded operators $L^\Delta_{\vartheta
 \vartheta''} : \mathcal{K}_{\vartheta''} \to \mathcal{K}_{\vartheta}$ we
define also a continuous linear operator $L^\Delta:\mathcal{K} \to
\mathcal{K} $, see (\ref{19c}). In view of this, we consider the
following two equations. First is
\begin{equation}
  \label{24c}
\frac{d}{dt} k_t = L^\Delta_\vartheta k_t, \qquad k_t|_{t=0} =
k_{\mu_0},
\end{equation}
considered as an equation in a given Banach space
$\mathcal{K}_{\vartheta}$. The second equation is (\ref{16}) with
$L^\Delta$ given in (\ref{15}) considered in the locally convex
space $\mathcal{K}$.
\begin{definition}
  \label{S1df}
By a solution of (\ref{24c}) on a time interval, $[0,T)$, $T\leq
+\infty$, we mean a continuous map $[0,T)\ni t \mapsto k_t \in
\mathcal{D} (L^\Delta_\vartheta)$ such that the map $[0,T)\ni t
\mapsto d k_t / dt\in \mathcal{K}_\vartheta$ is also continuous and
both equalities in (\ref{24c}) are satisfied. Likewise, a
continuously differentiable map $[0,T)\ni t \mapsto k_t \in
\mathcal{K}$ is said to be a solution of (\ref{16}) in $\mathcal{K}$
if both equalities therein are satisfied for all $t$. Such a
solution is called global if $T=+\infty$.
\end{definition}
\begin{remark}
  \label{D1rk}
The map $[0,T)\ni t\mapsto k_t \in \mathcal{K}$ is a solution of
(\ref{16}) if and only if, for each $t \in [0, T)$, there exists
$\vartheta''\in \mathds{R}$ such that $k_t\in
\mathcal{K}_{\vartheta''}$ and, for each $\vartheta > \vartheta''$,
the map $t\mapsto k_t$ is continuously differentiable at $t$ in
$\mathcal{K}_\vartheta$ and $d k_t/ dt = L^\Delta_\vartheta k_t =
L^\Delta_{\vartheta \vartheta''} k_t$.
\end{remark}
The main result of this subsection is contained in the following
statement.
\begin{theorem}
  \label{1tm}
Assume that (\ref{10}) and (\ref{11}) hold. Then for each $\mu_0 \in
\mathcal{P}_{\rm exp}(\Gamma^2)$, the problem (\ref{16}) with
(\ref{15}) with $k_0 = k_{\mu_0}$  has a unique global solution $k_t
\in \mathcal{K}^\star\subset \mathcal{K}$ which has the property
$k_t(\emptyset, \emptyset) = 1$. Therefore, for each $t\geq 0$ there
exists a unique state $\mu_t\in \mathcal{P}_{\rm exp}(\Gamma^2)$
such that $k_t = k_{\mu_t}$. Moreover, let $k_0$ and $C>0$ be such
that $k_0(\eta) \leq C^{|\eta_0|+ |\eta_1|}$ for $\lambda$-almost
all $\eta\in \Gamma_0^2$, see (\ref{9h}). Then the mentioned
solution satisfies
\begin{equation}
  \label{24d}
 \forall t\geq 0 \qquad 0\leq k_t (\eta) \leq C^{|\eta_0|+ |\eta_1|} \exp\left\{t\left(
 \alpha_0 |\eta_0| + \alpha_1 |\eta_1| \right) \right\}.
\end{equation}
\end{theorem}

\subsection{The mesoscopic level}

\label{SMe}

As is commonly recognized, see \cite[Chapter 8]{BL} and
\cite{Presutti}, the comprehensive understanding of the behavior of
an infinite interacting particle system requires its multi-scale
analysis. In the approach which we follow, see \cite{BKKK} (jump
dynamics) and \cite{FKKO} (two-component system), passing from the
micro- to the mesoscopic levels amounts to considering the system at
different spatial scales parameterized by $\varepsilon \in (0,1]$ in
such a way that $\varepsilon =1$ corresponds to the micro-level,
whereas the limit $\varepsilon \to 0$ yields the meso-level
description in which  the corpuscular structure disappears and the
system turns into a (two-component) medium characterized by a
density function $\varrho = (\varrho_0, \varrho_1)$, $\varrho_i :
\mathds{R}^d \to [0, +\infty)$, $i=0,1$. Then the evolution
$\varrho_0 \mapsto \varrho_t$, obtained from a {\it kinetic}
equation, approximates (in the mean-field sense) the evolution of
the system's states as it may be seen from the mesoscopic level.

\subsubsection{The kinetic equation}
Keeping in mind that the Poissonian state $\pi_\varrho$ is
completely characterized by the density $\varrho$, see (\ref{1fa})
and (\ref{1fb}), we introduce the following notion, cf. \cite[page
1046]{BKKK} and \cite[page 70]{FKKO}.
\begin{definition}
  \label{Ja1df}
A state $\mu\in \mathcal{P}_{\rm exp}(\Gamma^2)$ is said to be
Poisson-approximable if: (i)  there exist  $\vartheta\in \mathds{R}$
and $\varrho = (\varrho_0, \varrho_1)$, $\varrho_i \in L^\infty(
\mathds{R}^d \to \mathds{R})$, $\varrho_i \geq 0$, $i=0,1$, such
that both $k_\mu$ and $k_{\pi_\varrho}$ lie in
$\mathcal{K}_\vartheta$; (ii) for each $\varepsilon \in (0,1]$,
there exists $q_\varepsilon \in \mathcal{K}_\vartheta$ such that
$q_1 = k_\mu$ and $\|q_\varepsilon - k_{\pi_\varrho}\|_\vartheta \to
0$ as $\varepsilon \to 0$.
\end{definition}
Our aim is to show that the evolution $\mu_0 \mapsto \mu_t$ obtained
in Theorem \ref{1tm} preserves the property just defined relative to
the time dependent density $\varrho_t = (\varrho_{0,t},
\varrho_{1,t})$, obtained from the following system of kinetic
equations
\begin{eqnarray}
  \label{29}
  \left\{\begin{array}{ll}
\frac{d}{dt} \varrho_{0,t}  =(a_0 \ast \varrho_{0,t}) \exp\left(
- (\phi_0 \ast \varrho_{1,t})\right) \\[.2cm] \hskip5cm   - \varrho_{0,t} \left( a_0
\ast
\exp\left( - (\phi_0 \ast \varrho_{1,t})\right) \right),  \\[.5cm]
\frac{d}{dt} \varrho_{1,t} = (a_1 \ast \varrho_{1,t}) \exp\left( -
(\phi_1 \ast \varrho_{0,t})\right)\\[.2cm] \hskip5cm - \varrho_{1,t} \left( a_1 \ast
\exp\left( - (\phi_1 \ast \varrho_{0,t})\right) \right), \end{array}
\right. \qquad \quad
\end{eqnarray}
where $\ast$ denotes convolution; e.g.,
\begin{equation*}
  (a_i \ast \varrho_{i,t}) (x) = \int_{\mathds{R}^d} a_i (x-y)
  \varrho_{i,t} (y) d y, \qquad i=0,1.
\end{equation*}
\begin{definition}
  \label{Sksf}
By the global solution of the system of kinetic equations
(\ref{29}), subject to an initial condition,  we understand a
continuously differentiable map
\begin{equation}
  \label{40}
[0, +\infty) \ni t \mapsto (\varrho_{0,t}, \varrho_{1,t}) \in
L^\infty (\mathds{R}^d \to \mathds{R}^2)
\end{equation}
such that both equalities in (\ref{29}) hold. This solution is
called positive if $\varrho_{i,t}(x) \geq 0$, $i=0,1$, for all
$t\geq 0$ and Lebesgue-almost all $x\in \mathds{R}^d$. By the
positive solution of (\ref{29}) on the time interval $[0,T]$,
$0<T<\infty$, we mean the corresponding restriction of this map.
\end{definition}
Let $\|\cdot \|_{L^\infty}$ stand for the norm in $L^\infty
(\mathds{R}^d \to \mathds{R})$. In Theorem \ref{Ktm}, the space
$L^\infty (\mathds{R}^d \to \mathds{R}^2)$ is equipped with the norm
\begin{equation}
  \label{l-norm}
\| \varrho \|_{\infty} = \max_{i=0,1}\|\varrho_{i}\|_{L^\infty}.
\end{equation}
\begin{theorem}
  \label{Ktm}
For each positive $\varrho_0=(\varrho_{0,0}, \varrho_{1,0})\in
L^\infty (\mathds{R}^d \to \mathds{R}^2)$, the system of kinetic
equations (\ref{29}) with the initial condition $(\varrho_{0,t},
\varrho_{1,t})|_{t=0}=(\varrho_{0,0}, \varrho_{1,0})$  has a unique
positive global solution such that
\begin{equation}
  \label{40a}
\forall t \geq 0 \qquad  \varrho_{i,t} (x) \leq
\|\varrho_{i,0}\|_{L^\infty} \exp (\alpha_i t ), \quad i=0,1,
\end{equation}
where $\alpha_i$ are defined in (\ref{10}).
\end{theorem}
The relationship between the micro- and mesoscopic descriptions is
established by the following statement.
\begin{theorem}
  \label{4tm}
Let (\ref{33}) hold and $k_t$ and $\varrho_t$ be the solutions
described in Theorems \ref{1tm} and \ref{Ktm}, respectively. Assume
also that the initial state $\mu_0$ is Poisson-approximable by
$\pi_{\varrho_0}$, see Definition \ref{Ja1df}. That is, there exist
$\vartheta_*\in \mathds{R}$ and $q_{0, \varepsilon}$, $\varepsilon
\in (0,1]$, such that $k_{\mu_0} = q_{0,1}$ and $\|q_{0,\varepsilon}
- k_{\pi_{\varrho_0}}\|_{\vartheta_*} \to 0$ as $\varepsilon \to 0$.
Then there exist $\vartheta > \vartheta_*$ and $T>0$ such that
\begin{equation}
  \label{39}
 \lim_{\varepsilon \to 0} \sup_{t\in [0,T]} \|q_{t,\varepsilon} -
k_{\pi_{\varrho_t}}\|_{\vartheta} =0.
\end{equation}
\end{theorem}
Theorems \ref{Ktm} and \ref{4tm} are proved in Section \ref{S6}
below.

\subsubsection{The stationary solutions}
Stationary solutions $\varrho_{i,t} = \varrho_i$, $t\geq 0$, of the
system in (\ref{29}) are supposed to solve the following system of
equations
\begin{eqnarray}
  \label{Ch}
 \left\{ \begin{array}{ll}
(a_0 \ast \varrho_{0}) \exp\left( - (\phi_0 \ast \varrho_{1})\right)
 =  \varrho_{0} \left( a_0 \ast
\exp\left( - (\phi_0 \ast \varrho_{1})\right) \right), \\[.2cm]
(a_1 \ast \varrho_{1}) \exp\left( - (\phi_1 \ast \varrho_{0})\right)
=  \varrho_{1} \left( a_1 \ast \exp\left( - (\phi_1 \ast
\varrho_{0})\right) \right).  \end{array} \right.
\end{eqnarray}
It might be instructive to rewrite it in the form
\begin{equation}
  \label{CH1}
  \left\{ \begin{array}{ll} \psi_0 (x) = \int_{\mathds{R}^d} \tilde{a}_0
  (x,y) \psi_0 (y) dy, \\[.3cm]
\psi_1 (x) = \int_{\mathds{R}^d} \tilde{a}_1
  (x,y) \psi_1 (y) dy,
  \end{array}\right.
\end{equation}
where
\begin{eqnarray*}
\tilde{a}_0 (x,y) & := & \frac{a_0 (x-y) \exp\left( - (\phi_0 \ast
\varrho_{1})(y)\right)}{\int_{\mathds{R}^d}a_0 (x-y) \exp\left( -
(\phi_0 \ast \varrho_{1})(y)\right)d y}, \\[.2cm]
\tilde{a}_1 (x,y) & := & \frac{a_1 (x-y) \exp\left( - (\phi_1 \ast
\varrho_{0})(y)\right)}{\int_{\mathds{R}^d}a_1 (x-y) \exp\left( -
(\phi_1 \ast \varrho_{0})(y)\right)d y},
\end{eqnarray*}
and
\begin{equation}
  \label{CH}
 \psi_0 := \varrho_0 \exp\left(\phi_0 \ast
\varrho_{1}\right), \qquad  \psi_1 := \varrho_0 \exp\left( \phi_1
\ast \varrho_{0}\right).
\end{equation}
For each $\widetilde{C}_i> 0$, $i=0,1$, the system in (\ref{CH1})
has constant solutions $\psi_i \equiv \widetilde{C}_i$.  Then the
corresponding $\varrho_i$ are to be found from
\begin{equation}
  \label{Ch1}
\left\{\begin{array}{l} \varrho_0 = \widetilde{C}_0 \exp\left( -
(\phi_0 \ast \varrho_{1})\right), \\[.2cm] \varrho_1 = \widetilde{C}_1 \exp\left( -
(\phi_1 \ast \varrho_{0})\right). \end{array} \right.
\end{equation}
The solutions of (\ref{Ch1}) may be called \emph{birth-and-death}
solutions since they solve the corresponding equation of the
birth-and-death version of the Widom-Rowlinson dynamics with
specific values of $\widetilde{C}_i$, expressed in terms of the
model parameters, see \cite[eq. (4.13)]{FKKO}. The translation
invariant (i.e., constant) solution of (\ref{Ch1}) is
$\varrho_i\equiv C_i$, $i=0,1$, with $C_i$ satisfying, cf
(\ref{CH}),
\begin{equation}
  \label{Nc}
\widetilde{C}_0 = C_0 \exp\left( \langle \phi_0 \rangle C_1 \right),
\qquad \widetilde{C}_1 = C_1 \exp\left( \langle \phi_1 \rangle C_0
\right).
\end{equation}
For given $\widetilde{C}_0, \widetilde{C}_1>0$, let
$\mathcal{S}(\widetilde{C}_0, \widetilde{C}_1)$ be the set of all
positive $(\varrho_0,\varrho_1)\in L^\infty (\mathds{R}^d \to
\mathds{R}^2)$ that satisfy (\ref{Ch1}). Let also
$\mathcal{S}_c(\widetilde{C}_0, \widetilde{C}_1)$ be the subset of
$\mathcal{S}(\widetilde{C}_0, \widetilde{C}_1)$ consisting of
constant solutions $\varrho_i \equiv C_i$, $i=0,1$, with $C_i$
satisfying (\ref{Nc}). The symmetric case of (\ref{Nc}) with
specific values of $\widetilde{C}_i$  (as mentioned above) was
studied in \cite[Section 5]{FKKO}. Namely, for $ \langle \phi_1
\rangle \widetilde{C}_0 = \langle \phi_0 \rangle \widetilde{C}_1
=:a$, the set $\mathcal{S}_c(\widetilde{C}_0, \widetilde{C}_1)$ is a
singleton $\{C_0,C_1\}$ whenever $a\leq e$. Here
\begin{equation}
  \label{CH2}
C_0 = x_0/\langle \phi_1 \rangle , \qquad C_1 =x_0/ \langle
\phi_0\rangle ,
\end{equation}
with some $x_0 \in (0,1)$.  This solution is a stable node for
$a<e$. For $a>e$, there exist three solutions: (a) $ C_0 =
x_1/\langle \phi_1 \rangle$, $ C_1 = x_3/\langle \phi_0 \rangle$;
(b) $ C_0 = x_3/\langle \phi_1 \rangle$, $ C_1 = x_1/\langle \phi_0
\rangle$; (c) $ C_0 = x_2/\langle \phi_1 \rangle$,  $ C_1 =
x_2/\langle \phi_0 \rangle$. The first two solutions are stable
nodes and $x_3>1$. The stability means the existence of a small
neighborhood in $\mathcal{S}_c(\widetilde{C}_0, \widetilde{C}_1)$ of
the mentioned solution, which does not contain any other solution.

Let us now turn to the study of the stability of the constant
solutions of (\ref{Ch1}) with respect to perturbations $\varrho_i =
C_i + \epsilon_i$, $i=0,1$. By (\ref{Ch1}) and (\ref{Nc}) we
conclude that the perturbations ought to satisfy
\begin{equation}
  \label{Ch2}
\left\{\begin{array}{l} \epsilon_0 = C_0\left[ \exp\left\{-\left(
\phi_0 \ast  \epsilon_{1} \right) \right\} -1 \right], \\[.2cm] \epsilon_1 = C_1 \left[\exp\left\{-\left(
\phi_1 \ast \epsilon_{0} \right) \right\} -1 \right].
\end{array} \right.
\end{equation}
\begin{theorem}
  \label{K1tm}
The solution $\varrho_i\equiv  C_i$, $i=0,1$, of the system of
equations in (\ref{Ch}) is locally stable in
$\mathcal{S}(\widetilde{C}_0, \widetilde{C}_1)$, with
$\widetilde{C}_i$ and $C_i$ satisfying (\ref{Nc}), whenever the
following holds, cf. (\ref{CH2}),
\begin{equation}
  \label{Ch3}
C_0 C_1 \langle \phi_0 \rangle \langle \phi_1 \rangle < 1.
\end{equation}
This means that  there exists $\delta
>0$ such that $\varrho_i\equiv C_i$, $i=0,1$, is the only solution in the set $K_\delta:= \mathcal{S}(\widetilde{C}_0,
\widetilde{C}_1) \cap \{\varrho: \|\varrho-C\|_{\infty} < \delta\}$,
cf. (\ref{l-norm}).
\end{theorem}
\begin{proof}
Assume that $\|\epsilon_0\|_{L^\infty} >0$.  By means of the
inequality $|e^{-\alpha}-1| \leq |\alpha|e^{|\alpha|}$ we get from
(\ref{Ch2})
\[
\|\epsilon_0\|_{L^\infty} \leq  C_0 C_1 \langle \phi_0 \rangle
\langle \phi_1 \rangle \exp\left[\delta\left( \langle \phi_0 \rangle
+ \langle \phi_1 \rangle\right)\right]\cdot
\|\epsilon_0\|_{L^\infty} < \|\epsilon_0\|_{L^\infty},
\]
holding for small enough $\delta$ in view of (\ref{Ch3}). This
contradicts the assumption, and hence yields $\epsilon_0 =0$. The
corresponding estimate for $\|\epsilon_1\|_{L^\infty}$ is obtained
analogously.
\end{proof}
Assume now that both $\epsilon_i$  satisfy $\epsilon_i \in L^\infty
(\mathds{R}^d\to \mathds{R})\cap  L^1 (\mathds{R}^d\to \mathds{R})$.
Then each solution of (\ref{Ch2}) is a fixed point of the nonlinear
map $\Phi:L^\infty (\mathds{R}^d\to \mathds{R}^2)\cap L^1
(\mathds{R}^d\to \mathds{R}^2) \to L^\infty (\mathds{R}^d\to
\mathds{R}^2)\cap L^1 (\mathds{R}^d\to \mathds{R}^2) $ defined by
the right-hand of (\ref{Ch2}). Note that this $\Phi$ takes values in
$L^\infty (\mathds{R}^d\to \mathds{R}^2)\cap L^1 (\mathds{R}^d\to
\mathds{R}^2) $ in view of (\ref{33}). The zero solution of
(\ref{Ch2}) gets unstable whenever there exist nonzero $\epsilon =
(\epsilon_0, \epsilon_1)$ in the kernel of $I- \Phi'$, where $\Phi'$
is the Fr{\'e}chet derivative of $\Phi$ at $\epsilon = (0,0)$. By
(\ref{Ch2}) we have
\begin{equation}
  \label{CHh}
\Phi' \epsilon:= \Phi' \left(\!\begin{array}{ll} \epsilon_0\\[.2cm] \epsilon_1
\end{array} \!
\right) = \left(\!\begin{array}{ll} - C_0(\phi_0 \ast\epsilon_1)\\[.2cm] - C_1(\phi_1 \ast \epsilon_0)  \end{array}
\!\right).
\end{equation}
Since $\Phi'$ contains convolutions, it can be partially
diagonalized by means of the Fourier transform
\[
\hat{\phi}_i (p)= \int_{\mathds{R}^d} \phi_i (x) \exp\left(i
(p,x)\right)d x, \qquad p\in \mathds{R}^d, \ \ i=0,1.
\]
Note that both $\hat{\phi}_i$ are uniformly continuous on
$\mathds{R}^d$ and satisfy $|\hat{\phi}_i(p)| \leq \hat{\phi}_i(0) =
\langle \phi_i \rangle$, that  follows from their positivity.
Moreover, $|\hat{\phi}_i(p)| \to 0$ as $|p|\to +\infty$ (by the
Riemann-Lebesgue lemma). Note also that $\hat{\epsilon}_i$, $i=0,1$,
exist since $\epsilon_i$ are supposed to be integrable.
\begin{theorem}
  \label{K2tm}
Assume that the following holds, cf. (\ref{Ch3}),
\begin{equation}
  \label{N1c}
C_0 C_1 \langle \phi_0 \rangle \langle \phi_1 \rangle > 1.
\end{equation}
Then the constant solution $\varrho_i\equiv C_i$ of (\ref{Ch1}), and
hence of (\ref{Ch}), is unstable with respect to the perturbation
$\varrho_i = C_i + \epsilon_i$, $i=0,1$, with $\epsilon_i \in
L^\infty (\mathds{R}^d\to \mathds{R})\cap  L^1 (\mathds{R}^d\to
\mathds{R})$.
\end{theorem}
\begin{proof}
In view of the mentioned continuity of $\hat{\phi}_i$ and the
Riemann-Lebesgue lemma, the condition in (\ref{N1c}) implies the
existence of $p\in \mathds{R}^d\setminus \{0\}$ such that
\begin{equation}
  \label{N4c}
C_0 C_1 \hat{\phi}_0 (p) \hat{\phi}_1 (p) =1 .
\end{equation}
The instability in question takes place whenever the equation $\Phi'
\epsilon
 = \epsilon$, cf. (\ref{CHh}), has nonzero solutions in the
considered space. By means of the Fourier transform it can be turned
into
\begin{equation}
  \label{MNM}
\hat{\epsilon}_i (p) = C_0 C_1 \hat{\phi}_0 (p) \hat{\phi}_1 (p)
\hat{\epsilon}_i (p) , \qquad i=0,1,
\end{equation}
that has to hold for some $p\in \mathds{R}\setminus \{0\}$, which is
certainly the case in view of (\ref{N4c}).
\end{proof}
Given $C_i$, $i=0,1$, let $\epsilon =(\epsilon_0, \epsilon_1)$ solve
(\ref{Ch2}). Then $\varrho = (C_0+\epsilon_0, C_1+ \epsilon_1)$
solves (\ref{Ch1}) with $\widetilde{C}_i$ as in (\ref{Nc}) and hence
lies in $\mathcal{S}(\widetilde{C}_0, \widetilde{C}_1)$. Then
Theorem \ref{K2tm} describes the instability of the solution
$\varrho \equiv (C_0, C_1)$ in the latter set. For this reason, it
is independent of the jump kernels $a_i$. In order to study the
corresponding instability in the set of all solutions of (\ref{Ch}),
one has to rewrite (\ref{Ch}) in the form $\Psi(\varrho)=0$ and then
to show that the Fr{\'e}chet derivative $\Psi'$ of $\Psi$ at
$\varrho \equiv (C_0, C_1)$, defined as a bounded linear self-map of
$L^\infty (\mathds{R}^d\to \mathds{R})\cap  L^1 (\mathds{R}^d\to
\mathds{R})$, has nonzero $\epsilon$ in its kernel. By means of the
arguments used in the proof of Theorem \ref{K2tm} one readily
obtains that this is equivalent to, cf. (\ref{MNM}),
\[
\hat{\epsilon}_i (p) \left[1 - C_0 C_1 \hat{\phi}_0 (p) \hat{\phi}_1
(p)\right] \cdot \left[ \alpha_i - \hat{a}_i(p) \right] =0 , \qquad
i=0,1,
\]
that has to hold for some nonzero $p\in \mathds{R}^d$. Here
$\hat{a}_i(p)$, $i=0,1$, are the Fourier transforms of the jump
kernels, see (\ref{10}). Thus, if both these kernels are such that
$\hat{a}_i(p) < \hat{a}_i(0) = \alpha_i$ for all nonzero $p$, then
the latter condition turns into that in (\ref{MNM}).

\subsection{Comments}

\subsubsection{The microscopic description}
The only work on the Widom-Rowlinson dynamics of an infinite
particle system is that in \cite{FKKO} where a birth-and-death
(rather immigration-emigration) version was studied. In that
version, the particles of two types appear and disappear at random;
the appearance is subject to the repulsion from the particles of the
other type. The system's evolution was described by means of the
corresponding initial value problem for the Bogoliubov functional.
Namely, for $t< T$, where $T<\infty$ is expressed via the model
parameters, in \cite[Theorem 1]{FKKO} there was constructed the
evolution $B_{\mu_0} \mapsto B_t$, where $B_t: L^1(\mathds{R}^d \to
\mathds{R}^2)\to \mathds{R}$ is an exponential type entire function
and hence can be written down as, cf. (\ref{1e}),
\[
B_t (\theta) = \int_{\Gamma_0^2} k_t (\eta)E(\eta;\theta) \lambda (
d \eta).
\]
However, it was not shown that $B_t$ is the Bogoliubov functional,
i.e., that $k_t$ above is the correlation function, of some state
$\mu\in \mathcal{P}_{\rm exp} (\Gamma^2)$. In  the present work, for
the jump version of the Widom-Rowlinson model we show (Theorem
\ref{1tm}) that: (a) the evolution $k_{\mu_0} \mapsto k_t$, and
hence also $B_{\mu_0} \mapsto B_t$, can be continued to all $t>0$;
(b) for each $t>0$, $ B_t$ is the Bogoliubov functional of a unique
sub-Poissonian state $\mu_t$.

\subsubsection{The mesoscopic description}
In passing to the mesoscopic level of description, we use a scaling
procedure described in Section \ref{S4} below. It is equivalent to
the Lebowitz-Penrose scaling used in \cite{FKKO}, and also to the
Vlasov scaling used in \cite{BKKK,FKKK}. Our Theorem \ref{4tm} is
analogous to \cite[Theorem 2]{FKKO} proved for the birth-and-death
version. Note that the convergence in (\ref{39}) is uniform in $t$,
whereas in the mentioned statement of \cite{FKKO} the convergence is
point-wise.

Now we turn to the stationary solutions of (\ref{29}) which one
obtains from the system in (\ref{Ch}), or, equivalently, in
(\ref{CH1}). The latter may have nonconstant solutions $\psi_i$,
which then can be used to find the corresponding $\varrho_i$ from
(\ref{CH}). These solutions may depend on the jump kernels $a_i$.
The set of all solutions of (\ref{Ch}) contains the sets
$\mathcal{S}(\widetilde{C}_0, \widetilde{C}_1)$ for each pair
$\widetilde{C}_0$, $\widetilde{C}_1>0$. The corresponding solutions
$\varrho_i$ are independent of the jump kernels. Moreover,
$\mathcal{S}(\widetilde{C}_0, \widetilde{C}_1)$ is exactly the set
of solutions of the birth-and-death kinetic equation \cite[Eq.
(5.1)]{FKKO} corresponding to the specific values of
$\widetilde{C}_i$. Thus, our Theorems \ref{K1tm} and \ref{K2tm}
describe also the birth-and-death kinetic equation, which is an
extension of the study in \cite[Section 5]{FKKO}.

\section{The Rescaled Evolution}
\label{S4}

 In this section, we construct the evolution $q_{0,
\varepsilon} \mapsto q_{t, \varepsilon}$, $\varepsilon \in (0, 1]$,
which then will be used for: (a) obtaining the evolution stated in
Theorem \ref{1tm} in the form $k_t = q_{t,1}$; (b) proving Theorem
\ref{4tm}. To this end along with $L^\Delta$ defined in (\ref{15})
we will use
\begin{equation}
  \label{26}
  L^{\varepsilon,\Delta} = R^{-1}_\varepsilon  L^\Delta_{\varepsilon}
  R_\varepsilon, \qquad \varepsilon \in (0,1],
\end{equation}
where $L^\Delta_\varepsilon$ is obtained from $L^\Delta$ by
multiplying both $\phi_i$ by $\varepsilon$, and
\[
(R_\varepsilon q)(\eta_0, \eta_1) = \varepsilon^{-|\eta_0| -
|\eta_1|} q(\eta_0, \eta_1).
\]
We refer the reader to \cite{BKKK,FKKO} for more on deriving
operators as in (\ref{26}). Denote, cf. (\ref{12}),
\begin{equation}
  \label{31}
 \tau_{x,\varepsilon}^i (y) = \exp\left( -\varepsilon \phi_i (x-y)\right), \quad  t_{x,\varepsilon}^i
 (y) = \varepsilon^{-1} \left[\tau_{x,\varepsilon}^i (y) -1\right], \ \ i=0,1.
\end{equation}
Observe that
\begin{equation}
  \label{31a}
 \tau_{x,\varepsilon}^i (y) \to 1, \qquad t_{x,\varepsilon}^i
 (y) \to - \phi_i(x-y), \ \ {\rm as} \quad   \varepsilon \to 0.
\end{equation}
For $\varepsilon \in (0,1]$, let $Q^i_{y,\varepsilon}$ be as in
(\ref{13}) with $t_{x}^i$ replaced by $t_{x,\varepsilon}^i$ given in
(\ref{31}). Then the action of $L^{\varepsilon,\Delta}$ is  given by
the right-hand side of (\ref{15}) with both $Q^i_{y}$ replaced by
the corresponding $Q^i_{y,\varepsilon}$ and $\tau_{x}^i$ replaced by
$\tau_{x,\varepsilon}^i$. Note that, cf. (\ref{33}),
\begin{equation}
  \label{32}
\varepsilon^{-1}  \int_{\mathds{R}^d} \left( 1 - e^{-\varepsilon
\phi_i(x)}\right) d x \leq \langle \phi_i \rangle, \quad i=0,1.
\end{equation}
For each $\vartheta''\in \mathds{R}$, $k\in
\mathcal{K}_{\vartheta''}$, and $\varepsilon \in (0,1]$, by
(\ref{32}) both $Q^i_{y,\varepsilon}k$ satisfy the estimates as in
(\ref{21}) and (\ref{22}). Therefore, $L^{\varepsilon,\Delta} k$
satisfies (\ref{23}), which allows one to introduce the
corresponding linear operators $L^{\varepsilon,\Delta}_\vartheta :
\mathcal{D}(L^\Delta_\vartheta) \to \mathcal{K}_\vartheta$ and
$L^{\varepsilon,\Delta}_{\vartheta'\vartheta} :
\mathcal{K}_\vartheta \to \mathcal{K}_{\vartheta'}$, where
$\mathcal{D}(L^\Delta_\vartheta)$ is defined in (\ref{20}), see also
Corollary \ref{Gco} and (\ref{24a}). Thus, along with (\ref{24c}) we
will consider the problem
\begin{equation}
  \label{32a}
 \frac{d}{dt} q_{t,\varepsilon} = L^{\varepsilon,\Delta}_\vartheta
 q_{t,\varepsilon} , \qquad q_{t,\varepsilon}|_{t=0} =
 q_{0,\varepsilon}.
\end{equation}
Its solutions $q_{t,\varepsilon}\in
\mathcal{D}(L^\Delta_\vartheta)\subset \mathcal{K}_{\vartheta}$ are
defined analogously as in Definition \ref{S1df}.

For $\vartheta, \vartheta'\in \mathds{R}$ such that $\vartheta <
\vartheta'$, we set, cf. (\ref{24}),
\begin{gather}
  \label{21a}
T (\vartheta',\vartheta) = \frac{\vartheta' -
\vartheta}{4\alpha}\exp\left( - c
e^{\vartheta'}\right), \\[.2cm]
\alpha = \max_{i=0,1} \alpha_i, \qquad c = \max_{i=0,1} \langle
\phi_i \rangle.
 \nonumber
\end{gather}
For a fixed $\vartheta'\in \mathds{R}$, $T (\vartheta', \vartheta))$
can be made as big as one wants by taking small enough $\vartheta$.
However, if $\vartheta$ is fixed, then
\begin{equation}
  \label{N1}
 \sup_{\vartheta' > \vartheta} T(\vartheta', \vartheta) = \frac{\delta (\vartheta)}{4 \alpha} \exp\left( -
\frac{1}{\delta (\vartheta)}\right) =: \tau(\vartheta) < \infty,
\end{equation}
where $\delta(\vartheta)$ is the unique positive solution of the
equation
\begin{equation}
  \label{N2}
\delta e^\delta = \exp\left(- \vartheta - \log c \right).
\end{equation}
\begin{remark}
  \label{JJrk}
The supremum in (\ref{N1}) is attained at
\begin{equation*}
 \vartheta' = \vartheta + \delta (\vartheta).
\end{equation*}
Note also that $\delta (\vartheta) \to 0$, and hence
$\tau(\vartheta) \to 0$, as $\vartheta \to +\infty$.
\end{remark}
\begin{proposition}
  \label{3tm}
For arbitrary  $ \vartheta \in \mathds{R}$ and $\varepsilon \in
(0,1]$, the problem in (\ref{32a}) with $q_{0,\varepsilon}\in
\mathcal{K}_{\vartheta}$ has a unique solution $q_{t,\varepsilon}\in
\mathcal{K}_{\vartheta+ \delta(\vartheta)}$ on the time interval
$[0, \tau(\vartheta))$.
\end{proposition}
\begin{proof}
Take $T <  \tau(\vartheta)$ and then pick $\vartheta' \in
(\vartheta, \vartheta + \delta(\vartheta))$ such that $ T<
T(\vartheta', \vartheta)$. Our aim is to construct the family
\begin{equation}
  \label{40c}
S^\varepsilon_{\vartheta'\vartheta} (t) \in
\mathcal{L}(\mathcal{K}_{\vartheta}, \mathcal{K}_{\vartheta'}),
\qquad t\in [0, T ( \vartheta', \vartheta)),
\end{equation}
defined by the series
\begin{equation}
  \label{40b}
S^\varepsilon_{\vartheta'\vartheta} (t) = \sum_{n=0}^\infty
\frac{t^n}{n!} \left( L^{\varepsilon,
\Delta}\right)^n_{\vartheta'\vartheta}.
\end{equation}
In (\ref{40c}), $\mathcal{L}(\mathcal{K}_{\vartheta},
\mathcal{K}_{\vartheta'})$ stands for the Banach space of bounded
linear operators acting from $\mathcal{K}_{\vartheta}$ to
$\mathcal{K}_{\vartheta'}$ equipped with the corresponding operator
norm. In (\ref{40b}), $\left( L^{\varepsilon,
\Delta}\right)^0_{\vartheta'\vartheta}$ is the embedding operator
and
\begin{equation}
  \label{40d}
\left( L^{\varepsilon, \Delta}\right)^n_{\vartheta'\vartheta} :=
\prod_{l=1}^n L^{\varepsilon, \Delta}_{\vartheta_l \vartheta_{l-1}},
\quad \vartheta_l = \vartheta + l(\vartheta'- \vartheta)/n,
\end{equation}
for $n\in \mathds{N}$. Now we take into account that $\vartheta_l -
\vartheta_{l-1}= (\vartheta'- \vartheta)/n$ and that
$L^{\varepsilon, \Delta}$ satisfies (\ref{24}) for all $\varepsilon
\in (0,1]$. This yields the following estimate
\begin{eqnarray}
  \label{DA}
\|L^{\varepsilon, \Delta}_{\vartheta_l \vartheta_{l-1}}\| & \leq &
\left(\frac{n}{e}\right)(\vartheta' - \vartheta)\left\{ 2 \alpha_0
\exp\left( \langle \phi_0 \rangle e^{\vartheta'}\right) + 2 \alpha_1
\exp\left( \langle \phi_1 \rangle
e^{\vartheta'}\right)\right\}^{-1} \nonumber \\[.2cm]
& \leq & n \big{/}e T (\vartheta', \vartheta),
\end{eqnarray}
see (\ref{24}) and (\ref{21a}). Next we apply (\ref{DA}) in
(\ref{40d}) and conclude that the series in (\ref{40b}) converges in
the operator norm, uniformly on $[0,T]$, to the operator-valued
function $[0,T] \ni t \mapsto S^\varepsilon_{\vartheta'\vartheta}
(t) \in \mathcal{L}(\mathcal{K}_{\vartheta},
\mathcal{K}_{\vartheta'})$ such that
\begin{equation}
  \label{51}
\forall t\in [0,T]\qquad \|S^\varepsilon_{\vartheta'\vartheta} (t)
\| \leq \frac{T (\vartheta', \vartheta)}{T (\vartheta', \vartheta) -
t}.
\end{equation}
Likewise, for $\vartheta'' \in (\vartheta' , \vartheta + \delta
(\vartheta)]$, we get
\begin{eqnarray}
  \label{52}
 \frac{d}{dt} S^\varepsilon_{\vartheta''\vartheta} (t) & = &
 \sum_{n=0}^\infty \frac{t^n}{n!} \left( L^{\varepsilon,
\Delta}\right)^{n+1}_{\vartheta''\vartheta}\\[.2cm] & = & \sum_{n=0}^\infty
\frac{t^n}{n!} L^{\varepsilon, \Delta}_{\vartheta''
 \vartheta'} \left( L^{\varepsilon,
\Delta}\right)^n_{\vartheta'\vartheta} = L^{\varepsilon,
\Delta}_{\vartheta''
 \vartheta'} S^\varepsilon_{\vartheta'\vartheta} (t), \ \quad  t\in [0,T] \nonumber
\end{eqnarray}
Then
\begin{equation}
  \label{53}
 q_{t,\varepsilon} = S^\varepsilon_{\vartheta'\vartheta} (t) q_{0, \varepsilon} \in
\mathcal{K}_{\vartheta'} \subset
\mathcal{D}(L^{\varepsilon,\Delta}_{\vartheta''}),
\end{equation}
see Lemma \ref{1lm}, is a solution of (\ref{32a}) on the time
interval $[0, \tau(\vartheta))$ since $T< \tau(\vartheta)$ has been
taken in an arbitrary way.

Let us prove that the solution given in (\ref{53}) is unique. In
view of the linearity, to this end it is enough to show that the
problem in (\ref{32a}) with the zero initial condition has a unique
solution. Assume that $v_t\in
\mathcal{D}(L^{\varepsilon,\Delta}_{\vartheta +\delta(\vartheta)})$
is one of the solutions. Then $v_t$ lies in
$\mathcal{K}_{\vartheta''}$ for each $\vartheta'' > \vartheta +
\delta(\vartheta)$, see (\ref{19}). Fix any such $\vartheta''$ and
then take $t < \tau(\vartheta)$ such that $t< T (\vartheta'',
\vartheta+ \delta(\vartheta))$. Then, cf. (\ref{24a}),
\begin{eqnarray*}
v_t & = & \int_0^t L^{\varepsilon, \Delta}_{\vartheta''
 \bar{\vartheta}} v_s d s \\[.2cm]
& = & \int_0^t \int_0^{t_1} \cdots \int_0^{t_{n-1}} \left(
L^{\varepsilon, \Delta}\right)^n_{\vartheta''\bar{\vartheta}}
v_{t_n} d t_n \cdots d t_1, \nonumber
\end{eqnarray*}
where $\bar{\vartheta} := \vartheta + \delta(\vartheta)$ and $n\in
\mathds{N}$ is an arbitrary number. Similarly as above we get from
the latter
\[
\|v_t \|_{\vartheta''} \leq \frac{t^n}{n!} \left(\frac{n}{e T
(\vartheta'',\bar{ \vartheta})}\right)^n \sup_{s\in
[0,t]}\|v_s\|_{\bar{\vartheta}}.
\]
Since $n$ is an arbitrary number, this yields $v_s =0$ for all $s\in
[0,t]$. The extension of this result to all $t <\tau (\vartheta)$
can be done by repeating this procedure due times.
\end{proof}
\begin{remark}
  \label{Decrk}
Similarly as in obtaining (\ref{52}) we have that, for each
$\varepsilon \in (0,1]$ and all $\vartheta_0, \vartheta_1,
\vartheta_2 \in \mathds{R}$ such that $\vartheta_0 < \vartheta_1 <
\vartheta_2$, the following holds
\begin{eqnarray}
  \label{D1}
& & \quad S^\varepsilon_{\vartheta_2 \vartheta_0} (t+s) =
S^\varepsilon_{\vartheta_2 \vartheta_1} (t)
S^\varepsilon_{\vartheta_1
\vartheta_0} (s), \\[.2cm] & &  t\in [0, T(\vartheta_2, \vartheta_1)), \quad s \in [0,
T(\vartheta_1, \vartheta_0)). \nonumber
\end{eqnarray}
\end{remark}

\section{The Proof of Theorem \ref{1tm}}
\label{SPr}

With the help of Proposition \ref{3tm} we have already obtained the
unique solution of (\ref{24c}) in the form \begin{equation}
  \label{JN1}
k_t = S^1_{\vartheta \vartheta_0} k_{\mu_0}, \qquad t< \tau
(\vartheta_0),
\end{equation}
where $k_{\mu_0} \in \mathcal{K}_{\vartheta_0}$ and $\vartheta \in(
\vartheta_0 , \vartheta_0 + \delta(\vartheta_0))$ is taken such that
$t < T(\vartheta_0 + \delta(\vartheta_0),\vartheta)$. To prove
Theorem \ref{1tm} we first show in Lemma \ref{Id1lm}  that $k_t$
lies in the cone (\ref{19a}) and hence is a correlation function of
a unique state $\mu_t$. Then in Lemma \ref{A1lm} we construct an
auxiliary evolution $u_0 \mapsto u_t$, with which we compare the
evolution $k_{\mu_0} \mapsto k_t$ defined in (\ref{JN1}). Finally,
we construct its extension to all $t>0$ as stated in the theorem.

\subsection{The identification lemma}

Our aim now is to show that the solution of (\ref{24c}) given in
(\ref{JN1}) has the property $k_t \in
\mathcal{K}^\star_{\vartheta}$, see (\ref{19a}).
\begin{lemma}
  \label{Id1lm}
Let $\vartheta$ and $\vartheta^*$ be as in Corollary \ref{Gco}. Then
for each $t\in [0, T(\vartheta, \vartheta^*))$, the operator defined
in (\ref{40b}) has the property
\begin{equation}
  \label{63}
  S^1_{\vartheta \vartheta^*}(t):\mathcal{K}^\star_{\vartheta^*} \to
  \mathcal{K}^\star_{\vartheta}.
\end{equation}
\end{lemma}
\begin{proof}
We follow the line of arguments used in the proof of Theorem 3.8 of
\cite{BKKK}, see also \cite[Lemma 4.8]{KK}. Let $\mu_0\in
\mathcal{P}_{\rm exp} (\Gamma^2)$ be such that $k_{\mu_0} \in
\mathcal{K}_{\vartheta^*}^\star$, see Proposition \ref{Gpn}. For
$\Lambda=(\Lambda_0,\Lambda_1)$, $\Lambda_i\in \mathcal{B}_{\rm
b}(\mathds{R}^d)$, $i=0,1$, let $\mu^\Lambda_0$ and
$R^\Lambda_{\mu_0}$ be as in (\ref{9c}). For $N\in \mathds{N}$, we
then set
\begin{equation}
  \label{64}
R^{\Lambda,N}_0 (\eta) = R^\Lambda_{\mu_0} (\eta) I_N (\eta), \qquad
\eta \in \Gamma_0^2,
\end{equation}
where $I_N (\eta)=1$ whenever $\max_{i=0,1} |\eta_i| \leq N$ and
$I_N (\eta)=0$ otherwise. Set
\begin{eqnarray}
  \label{65}
\mathcal{R} & = & L^1 (\Gamma^2_0, d \lambda) , \quad
\mathcal{R}_\beta = L^1 (\Gamma^2_0, b_\beta d \lambda), \\[.2cm]
b_\beta (\eta) & := & \exp\bigg{(} \beta\left(  |\eta_0| +
|\eta_1|\right)\bigg{)}, \qquad \beta >0.\nonumber
\end{eqnarray}
Let $\|\cdot\|_{\mathcal{R}}$ and $\|\cdot\|_{\mathcal{R}_\beta}$ be
the norms of the spaces introduced in (\ref{65}) and $\mathcal{R}^+$
and $\mathcal{R}^+_\beta$ be the corresponding cones of positive
elements. For each $\beta>0$, $R^{\Lambda,N}_0$ defined in
(\ref{64}) lies in $\mathcal{R}_\beta^+ \subset \mathcal{R}^+$ and
is such that $\|R^{\Lambda,N}_0\|_{\mathcal{R}} \leq 1$. Similarly
as for the Kawasaki model, see \cite[Section 3.2]{BKKK}, it is
possible to show that $L^*$ related  by (\ref{1b}) to $L$ given in
(\ref{LF}) generates the evolution of states $\mu_0 \mapsto \mu_t$,
$t\geq 0$, whenever $\mu_0$ has the property $\mu_0(\Gamma^2_0)=1$,
which is the case for  $\mu^\Lambda_0$. Moreover, for each $t\geq
0$, the mentioned  $\mu_t$ is absolutely continuous with respect to
$\lambda$, and the equation for $R_t = d\mu_t / d \lambda$
corresponding to (\ref{1a}) can be written in the form
\begin{equation}
  \label{66}
\frac{d}{dt} R_t = L^\dagger R_t, \qquad R_t|_{t=0}=R_{\mu_0},
\end{equation}
where, cf. (\ref{15}), $L^\dagger$ is defined by the relation
$L^\dagger R = d(L^* \mu)/d\lambda$, and hence acts according to the
following formula
\begin{eqnarray}
  \label{66a}
( L^\dagger R)(\eta_0, \eta_1)& = & \sum_{y\in \eta_0}
\int_{\mathds{R}^d} a_0 (x-y) e(\tau^0_y;\eta_1) R(\eta_0\setminus y
\cup x, \eta_1) d x \qquad \\[.2cm]
& + & \sum_{y\in \eta_1} \int_{\mathds{R}^d} a_1 (x-y)
e(\tau^1_y;\eta_0) R(\eta_0, \eta_1\setminus y
\cup x) d x \nonumber \\[.2cm] & - & \Psi(\eta_0, \eta_1) R(\eta_0,
\eta_1), \nonumber \\[.2cm]
\Psi(\eta_0, \eta_1)& := & \sum_{x\in \eta_0} \int_{\mathds{R}^d}
a_0 (x-y) e(\tau^0_y;\eta_1) d y \nonumber \\[.2cm]
& & \qquad + \sum_{x\in \eta_1} \int_{\mathds{R}^d} a_1 (x-y)
e(\tau^1_y;\eta_0) d y. \nonumber
\end{eqnarray}
Like in \cite[Theorem 3.7]{BKKK}, one shows that $L^\dagger$
generates a stochastic $C_0$-semigroup, $S_R:= \{S_R(t)\}_{t\geq
0}$, on $\mathcal{R}$, which leaves invariant each
$\mathcal{R}_\beta$, $\beta >0$. Then the solution of (\ref{66}) is
$R_t = S_R (t) R_0$. For $R_0^{\Lambda,N}$ as in (\ref{64}), we then
set
\begin{equation}
  \label{67}
R_t^{\Lambda,N} (t) = S_R(t) R_0^{\Lambda,N}, \qquad t>0.
\end{equation}
Then $R_t^{\Lambda,N}\in \mathcal{R}_\beta^+ \subset \mathcal{R}^+$
and $\|R^{\Lambda,N}_t\|_{\mathcal{R}} \leq 1$. This yields that,
for each $G\in B_{\rm bs}^\star (\Gamma_0^2)$, see (\ref{9f}) and
(\ref{9g}), the following holds
\begin{equation}
  \label{68}
 \langle \! \langle KG , R^{\Lambda,N}_t \rangle\!\rangle \geq 0,
 \qquad t\geq 0.
\end{equation}
The integral in (\ref{68}) exists as $R_t^{\Lambda,N}\in
\mathcal{R}_\beta$ and $KG$ satisfies (\ref{9b}). Moreover, like in
(\ref{24}) and (\ref{56}), for each $\beta'$ such that $0< \beta' <
\beta$, we derive from (\ref{66a}) the following estimate
\begin{equation*}
\|L^\dagger R\|_{\mathcal{R}_{\beta'}} \leq \frac{4\alpha \|
R\|_{\mathcal{R}_{\beta}} }{e(\beta - \beta')}.
\end{equation*}
This allows us to define the corresponding bounded operators
$(L^\dagger)^n_{\beta'\beta} : \mathcal{R}_{\beta} \to
\mathcal{R}_{\beta'}$, $n\in \mathds{N}$, cf. (\ref{40d}) and
(\ref{59}), the norms of which satisfy
\begin{equation}
  \label{68b}
\|(L^\dagger)^n_{\beta'\beta} \| \leq n^n \left(
e\bar{T}(\beta,\beta')\right)^{-n}.
\end{equation}
On the other hand, we have that, cf. (\ref{9d}) and (\ref{64}),
\begin{eqnarray}
  \label{69}
k_0^{\Lambda,N} (\eta)& := & \int_{\Gamma_0^2} R^{\Lambda,N}_0
(\eta\cup\xi)\lambda(d\xi)\\[.2cm]&= & \int_{\Gamma_0^2} R^{\Lambda,N}_0
(\eta_0\cup\xi_0, \eta_1\cup\xi_1)(\lambda_0\otimes
\lambda_1)(d\xi_0, d \xi_1) \nonumber
\end{eqnarray}
is such that $k_0^{\Lambda,N}\in \mathcal{K}_{\vartheta^*}^\star$,
and hence we may get
\begin{equation}
  \label{70}
k_t^{\Lambda,N} = S^1_{\vartheta \vartheta^*}(t) k_0^{\Lambda,N},
\qquad t\in [0,T(\vartheta, \vartheta^*)),
\end{equation}
where $S^1_{\vartheta \vartheta^*}(t) = S^\varepsilon_{\vartheta
\vartheta^*}(t)|_{\varepsilon=1}$ is given in (\ref{40b}).
 Then the proof of (\ref{63}) consists in showing:
\begin{eqnarray}
  \label{71}
 &(i) & \quad \forall G\in B^\star_{\rm bs}(\Gamma^2_0) \qquad
 \langle \! \langle G, k^{\Lambda,N}_t \rangle \!\rangle \geq
 0;\\[.2cm]
 &(ii)& \quad  \langle \!\langle G, S^1_{\vartheta \vartheta^*} (t) k_0 \rangle \!\rangle =
 \lim_{\Lambda \to \mathds{R}^d \times \mathds{R}^d} \lim_{N\to +\infty} \langle \! \langle G, k^{\Lambda,N}_t \rangle \!\rangle .\nonumber
\end{eqnarray}
To prove claim {\it (i)} of (\ref{71}), for a given $G\in
B^\star_{\rm bs}(\Gamma^2_0)$ one sets
\begin{equation}
\label{72} \varphi_G (t) = \langle \! \langle KG , R^{\Lambda,N}_t
\rangle\!\rangle, \quad \ \  \psi_G (t) = \langle \! \langle G ,
k^{\Lambda,N}_t \rangle\!\rangle,
\end{equation}
where $\psi_G$ is defined for $t$ as in (\ref{70}). For a given
$t\in (0,T (\vartheta, \vartheta^*))$, we pick
$\vartheta'<\vartheta$ such that $t < T (\vartheta', \vartheta^*)$,
and hence $k_{s}^{\Lambda,N} \in \mathcal{K}_{\vartheta'}$ for $s\in
[0,t]$. Then the  direct calculation based on (\ref{52}) yields for
the $n$-th derivative
\[
\psi_G^{(n)} (t) = \langle \! \langle G, (L^\Delta)^n_{\vartheta
\vartheta'} k^{\Lambda,N}_t \rangle \!\rangle, \qquad n \in
\mathds{N}.
\]
As in obtaining (\ref{51}) we then get from the latter
\begin{equation}
  \label{73}
|\psi_G^{(n)} (t) | \leq A^n n^n C_{\vartheta'}(G) \sup_{s\in
[0,t]}\|k^{\Lambda,N}_s\|_{\vartheta'}.
\end{equation}
Here $A= 1/ e T(\vartheta, \vartheta')$ and
\[
C_{\vartheta'}(G)  = \int_{\Gamma_0^2} |G(\eta)|\exp\left( \vartheta
' |\eta_0| + \vartheta '|\eta_1| \right) \lambda(d\eta) <\infty,
\]
as $G\in B_{\rm bs}(\Gamma_0^2)$, see Definition \ref{Gdef}.
Likewise, from  (\ref{67}) we have
\[
\varphi^{(n)}_G (t) = \langle \! \langle KG,
(L^\dagger)^n_{\beta'\beta} R^{\Lambda,N}_t \rangle \! \rangle
\]
For the same $t$ as in (\ref{73}), by (\ref{68b}) we have from the
latter
\begin{equation}
  \label{73a}
|\varphi_G^{(n)} (s) | \leq \bar{A}^n n^n C_{\beta'}(G) \sup_{s\in
[0,t]}\|R^{\Lambda,N}_s\|_{\beta'}.
\end{equation}
Here $\bar{A}= 1/ e \bar{T}(\beta', \beta)$ and
\[
C_{\beta'}(G)  =  \esssup_{\eta \in \Gamma_0^2} |KG(\eta)| \exp
\left(-\beta' |\eta_0| - \beta'|\eta_1|\right) <\infty
\]
which holds in view of (\ref{9b}). By (\ref{15}), (\ref{66a}), and
(\ref{69}) it follows that
\[
(L^\Delta k^{\Lambda,N}_0)(\eta)=\int_{\Gamma^2_0}(L^\dagger
R_0^{\Lambda,N})(\eta \cup \xi) \lambda (d\xi),
\]
which then yields
\begin{equation}
  \label{74}
\forall n\in \mathds{N}_0 \qquad \varphi^{(n)}_G(0) =
\psi^{(n)}_G(0).
\end{equation}
By the Denjoy-Carleman theorem \cite{DC}, (\ref{73a}) and (\ref{73})
imply that both functions defined in (\ref{72}) are quasi-analytic
on $[0,t]$. Then (\ref{74}) implies
\begin{equation}
  \label{74a}
 \forall t\in [0,T(\vartheta, \vartheta^*)) \qquad \varphi_G(s) = \psi_G(s),
\end{equation}
which by (\ref{68}) yields the first  line in (\ref{71}). The
convergence claimed in {\it(ii)} of (\ref{71}) is proved in a
standard way, see Appendix in \cite{BKKK}.
\end{proof}
Note that (\ref{74a}) yields also that
\begin{equation}
  \label{74b}
 \forall s\in [0,T(\vartheta, \vartheta^*)) \qquad \langle
 \!\langle G,  q^{\Lambda,N}_t\rangle \!\rangle
= \langle
 \!\langle G,  k^{\Lambda,N}_t\rangle \!\rangle,
\end{equation}
where $G$ and  $k^{\Lambda,N}_t$ are as in (\ref{72}) and
\begin{equation}
  \label{74c}
q^{\Lambda,N}_t (\eta) := \int_{\Gamma_0^2} R^{\Lambda,N}_t (\eta
\cup \xi) \lambda (d\xi),
\end{equation}
cf. (\ref{69}).

\subsection{An auxiliary evolution}

The evolution which we construct now  will be used to extending the
solution $k_t$ given in (\ref{JN1}) to the global solution as stated
in Theorem \ref{1tm}. The construction employs the operator
\begin{eqnarray}
  \label{55}
(\bar{L} k) (\eta_0, \eta_1) & = & \sum_{y\in \eta_0}
\int_{\mathds{R}^d} a_0 (x-y) k(\eta_0 \setminus y \cup x, \eta_1) d
x \\[.2cm] & + & \sum_{y\in \eta_1}
\int_{\mathds{R}^d} a_1 (x-y) k(\eta_0, \eta_1 \setminus y \cup x) d
x  \nonumber
\end{eqnarray}
obtained from $L^\Delta$ given in (\ref{15}) by putting $\phi_i =0$,
$i=0,1$, and then dropping the second and fourth terms. Hence, like
in (\ref{24}) we get
\begin{equation}
  \label{56}
\| \bar{L}k\|_\vartheta \leq  \frac{4 \alpha \|k\|_{\vartheta''}}{
e(\vartheta - \vartheta'')},
\end{equation}
which allows us to introduce the operators $(\bar{L}_\vartheta,
\mathcal{D} (\bar{L}_\vartheta))$ and
$\bar{L}_{\vartheta\vartheta''}\in
\mathcal{L}(\mathcal{K}_{\vartheta''}, \mathcal{K}_{\vartheta})$
such that, cf. (\ref{24a}),
\begin{equation*}
\forall k \in \mathcal{\vartheta''} \qquad
\bar{L}_{\vartheta\vartheta''}k = \bar{L}_{\vartheta} k, \qquad
\vartheta'' < \vartheta.
\end{equation*}
Like above, we have that
\[
\mathcal{K}_{\vartheta''} \subset \mathcal{D}(\bar{L}_{\vartheta}):=
\{ k \in \mathcal{K}_{\vartheta} : \bar{L} k \in
\mathcal{K}_{\vartheta}\}, \qquad \vartheta'' < \vartheta  .
\]
Note that
\begin{equation}
  \label{58}
\bar{L}_{\vartheta\vartheta''} :\mathcal{K}_{\vartheta''}^+ \to
\mathcal{K}_{\vartheta}^+ , \qquad \vartheta'' < \vartheta,
\end{equation}
see (\ref{19b}).  For $n\in \mathds{N}$, we define
$(\bar{L})^n_{\vartheta'\vartheta}$ similarly as in (\ref{40d}) and
denote, cf. (\ref{21a}),
\begin{equation}
  \label{59}
\bar{T}(\vartheta', \vartheta) = (\vartheta' - \vartheta)/4 \alpha,
\qquad \vartheta < \vartheta' .
\end{equation}
Our aim is to study the operator valued function defined by the
series
\begin{equation}
  \label{61}
\bar{S}_{\vartheta'\vartheta} (t) = \sum_{n=0}^\infty \frac{t^n}{n!}
\left(\bar{L}\right)^n_{\vartheta'\vartheta}.
\end{equation}
\begin{lemma}
  \label{A1lm}
For each $\vartheta_0, \vartheta\in \mathds{R}$ such that
$\vartheta_0 < \vartheta$, the series in (\ref{61}) defines a
continuous function
\begin{equation}
  \label{DD}
[0, \bar{T}(\vartheta, \vartheta_0) ) \ni t \mapsto
\bar{S}_{\vartheta\vartheta_0} (t) \in \mathcal{L}(
\mathcal{K}_{\vartheta_0}, \mathcal{K}_{\vartheta}),
\end{equation}
which has the following properties:
\begin{itemize}
  \item[{\it (a)}] For $t$ as in (\ref{DD}), let $\vartheta''\in (\vartheta_0, \vartheta)$ be such that
 $t< \bar{T}(\vartheta'',
\vartheta_0)$. Then, cf. (\ref{52}),
\begin{equation}
  \label{D2}
\frac{d}{dt}\bar{S}_{\vartheta\vartheta_0}(t) = \bar{L}_{\vartheta
\vartheta''} \bar{S}_{\vartheta''\vartheta_0} (t).
\end{equation}
  \item[{\it (b)}]
The problem
\begin{equation}
  \label{60}
 \frac{d}{dt} u_t = \bar{L}_\vartheta u_t , \qquad u_t|_{t=0} = u_0
 \in \mathcal{K}^+_{\vartheta_0},
\end{equation}
has a unique solution $u_t \in \mathcal{K}^+_{\vartheta}$ on the
time interval $[0, \bar{T}(\vartheta, \vartheta_0))$ given by
\begin{equation}
  \label{62}
  u_t = \bar{S}_{\vartheta''\vartheta_0} (t)u_0,
\end{equation}
where, for a fixed $t \in [0, \bar{T}(\vartheta, \vartheta_0))$,
$\vartheta''$ is chosen to satisfy $t< \bar{T}(\vartheta'',
\vartheta_0)$.
\end{itemize}
\end{lemma}
\begin{proof}
Proceeding as in the proof of Proposition \ref{3tm}, by means of the
estimate in (\ref{56}) we prove the convergence of the series in
(\ref{61}). This allows also for proving (\ref{D2}), which yields
the existence of the solution of (\ref{60}) in the form given in
(\ref{62}). The uniqueness is proved analogously as in the case of
Proposition \ref{3tm}. The stated positivity of $u_t$ follows from
(\ref{61}) and (\ref{58}).
\end{proof}
\begin{corollary}
  \label{D1co}
For a given $C>0$, we let in (\ref{60}) and (\ref{62})  $\vartheta_0
=\log C$ and $u_0 (\eta) = C^{|\eta_0|+|\eta_1|}$. Then the unique
solution of (\ref{60}) is
\begin{equation}
  \label{85}
u_t (\eta) = C^{|\eta_0|+|\eta_1|} \exp\left\{ t(\alpha_0 |\eta_0|+
\alpha_1 |\eta_1|)\right\}.
\end{equation}
This solution can naturally be continued to all $t>0$ for which it
lies in $\mathcal{K}_{\vartheta(t)}$ with
\begin{equation}
  \label{De}
\vartheta(t) =  \log C + t \max_{i=0,1}\alpha_i.
\end{equation}
\end{corollary}
\begin{proof}
In view of the lack of interaction in (\ref{55}), the equations for
particular $u^{(n)}_t$ take the following (decoupled) form
\begin{eqnarray*}
& &  \frac{d}{dt} u_t^{(n)} (x_1, \dots ,x_{n_0}; y_1, \dots ,y_{n_1}) \\[.2cm]
& & \quad =
 \sum_{i=1}^{n_0} \int_{\mathds{R}^d} a_0 (x- x_i) u_t^{(n)} (x_1,
 \dots , x_{i-1}, x , x_{i+1}, \dots ,x_{n_0}; y_1, \dots ,y_{n_1}) d x
 \\[.2cm]
 & &  +
 \sum_{i=1}^{n_1} \int_{\mathds{R}^d} a_1 (y- y_i) u_t^{(n)} (x_1,
 \dots x_{n_0};  y_1, \dots y_{i-1}, y , y_{i+1}, \dots  y_{n_1}) d
 y, \  n\in  \mathds{N}^2, \qquad \qquad
\end{eqnarray*}
which for the initial translation invariant $u_0$ yields (\ref{85}).
\end{proof}

\subsection{The global solution}

As follows from Proposition \ref{3tm} and Lemma \ref{Id1lm}, the
unique solution of the problem (\ref{24c}) with $k_{0}\in
\mathcal{K}^\star_{\vartheta^*}$ lies in
$\mathcal{K}_\vartheta^\star$ for $t \in (0, T (\vartheta,
\vartheta^*))$. At the same time, for fixed $\vartheta^*$,  $T
(\vartheta, \vartheta^* )$ is bounded, see (\ref{N1}). This means
that the mentioned solution cannot be directly continued as stated
in Theorem \ref{1tm}. In this subsection,  by a comparison method we
prove that, for $t \in (0, T (\vartheta, \vartheta^*))$, $k_t$
satisfies (\ref{24d}) which is then used to get the continuation in
question, cf. Corollary \ref{D1co}. Recall that the operators
$Q_y^i$, $i=0,1$, were introduced in (\ref{13}) and the cone
$\mathcal{K}^+_\vartheta$ was defined in (\ref{19b}).
\begin{lemma}
  \label{J1lm}
For each $k_0\in \mathcal{K}_{\vartheta^*}^\star$ and $t \in (0, T
(\vartheta, \vartheta^*))$, $k_t := S^1_{\vartheta \vartheta^*}k_0$
has the property
\begin{equation}
  \label{75}
k_t - e(\tau^i_y;\cdot) Q^i_y k_t \in \mathcal{K}_\vartheta^+,
\qquad i=0,1,
\end{equation}
holding for Lebesgue-almost all $y\in \mathds{R}^d$.
\end{lemma}
\begin{proof}
Clearly, it is enough to show that (\ref{75}) holds for $i=0$. For a
fixed $y$, we denote
\[
v_{t,1} = k_t - Q^0_y k_t, \quad v_{t,2} =  [1- e(\tau^0_y;\cdot)]
Q^0_y k_t.
\]
The proof will be done if we show that, for all $G\in B_{\rm
bs}(\Gamma^2_0)$ such that $G(\eta) \geq 0$ for $\lambda$-almost all
$\eta\in \Gamma_0^2$, the following holds
\begin{equation}
  \label{76}
\langle \! \langle G, v_{t,j} \rangle \! \rangle \geq 0, \qquad
j=1,2.
\end{equation}
Let $\Lambda$, $N$, and $k_0^{\Lambda,N}$ be as in (\ref{69}), and
then $k_t^{\Lambda,N}$ be as in (\ref{70}). Next, let
$v_{t,j}^{\Lambda,N}$, $j=1,2$, be defined as above with $k_t$
replaced by $k_t^{\Lambda,N}$. By (\ref{74b}) and (\ref{74c}) we
then get
\begin{eqnarray}
  \label{78}
\langle \! \langle G, Q^0_y k^{\Lambda,N}_t \rangle \! \rangle & = &
\int_{\Gamma_0^2} \widetilde{G}(\eta) k^{\Lambda,N}_t (\eta)
\lambda(d \eta) \\[.2cm]
& = & \int_{\Gamma_0^2} \int_{\Gamma_0^2} \widetilde{G}(\eta)
R^{\Lambda,N}_t (\eta \cup \xi) \lambda (d\eta) \lambda (d\xi),
\nonumber
\end{eqnarray}
where
\begin{equation*}
  \widetilde{G}(\eta_0,\eta_1):=\sum_{\xi\subset \eta_1} e(t^0_y; \xi) G(\eta_0 , \eta_1 \setminus
\xi).
\end{equation*}
Furthermore, by (\ref{78}) we get
\begin{eqnarray}
  \label{78b}
& & \langle \! \langle G, Q^0_y k^{\Lambda,N}_t \rangle \! \rangle
=  \int_{\Gamma_0^2} G(\eta_0, \eta_1) \\[.2cm] & & \quad
 \times \int_{\Gamma_0^2} \left(
\int_{\Gamma_0} e(t^0_y;\zeta) R^{\Lambda,N}_t(\eta_0\cup \xi_0,
\eta_1 \cup \xi_1 \cup \zeta) \lambda_1 ( d \zeta) \right)
\lambda(d\eta) \lambda ( d \xi) \nonumber \\[.2cm]
& & \quad =  \int_{\Gamma_0^2} G(\eta_0, \eta_1) \int_{\Gamma_0^2}
\left( \sum_{\zeta\subset \xi_1} e(t^0_y;\zeta)\right)
R^{\Lambda,N}_t(\eta_0\cup \xi_0, \eta_1 \cup \xi_1) \lambda(d\eta)
\lambda ( d \xi).\nonumber \qquad
\end{eqnarray}
By (\ref{12}) we have that
\begin{equation*}
 \sum_{\zeta \subset \xi_1} e(t^0_y;\zeta) = e(\tau^0_y;\xi_1).
\end{equation*}
We apply this in the last line of (\ref{78b}) and obtain
\begin{eqnarray}
  \label{78c}
& & \langle \! \langle G, Q^0_y k^{\Lambda,N}_t \rangle \! \rangle
\\[.2cm] & & \quad =  \int_{\Gamma_0^2} G(\eta_0, \eta_1) \int_{\Gamma_0^2}
 e(\tau^0_y;\xi_1)
R^{\Lambda,N}_t(\eta_0\cup \xi_0, \eta_1 \cup \xi_1) \lambda(d\eta)
\lambda ( d \xi)\nonumber \qquad \\[.2cm]
& & \quad \leq \int_{\Gamma_0^2} G(\eta_0, \eta_1) \int_{\Gamma_0^2}
 R^{\Lambda,N}_t(\eta_0\cup \xi_0, \eta_1 \cup \xi_1) \lambda(d\eta)
\lambda ( d \xi)\nonumber \qquad \\[.2cm] & & \quad = \langle \! \langle G,  k^{\Lambda,N}_t \rangle \!
\rangle, \nonumber
\end{eqnarray}
which after the limiting transition as in (\ref{71}) yields
(\ref{76}) for $j=1$. For the same $G$, we set $\bar{G} =
e(\tau^0_y;\cdot) G$. Then by (\ref{12}) and the second line in
(\ref{78c}) we get
\begin{equation*}
\langle \! \langle \bar{G}, Q^0_y k^{\Lambda,N}_t \rangle \! \rangle
\leq \langle \! \langle G, Q^0_y k^{\Lambda,N}_t \rangle \! \rangle,
\end{equation*}
which after the limiting transition as in (\ref{71}) yields
(\ref{76}) for $j=2$.
\end{proof}
\begin{lemma}
  \label{J2lm}
Let $C>0$ be such that the initial condition in (\ref{24c})
satisfies $k_{\mu_0}(\eta) =k_0 (\eta) \leq C^{|\eta_0|+|\eta_1|}$.
Then for all $t< T (\vartheta, \vartheta^*)$ with $\vartheta^* =
\log C$ and any $\vartheta> \vartheta^*$, the unique solution of
(\ref{24c}) given by the formula
\begin{equation}
  \label{Sol}
k_t = S^1_{\vartheta \vartheta^*} (t) k_0
\end{equation}
satisfies (\ref{24d}) for $\lambda$-almost all $\eta\in \Gamma_0^2$.
\end{lemma}
\begin{proof}
Take any $\vartheta > \vartheta^*$ and fix $t< T (\vartheta,
\vartheta^*)$; then pick $\vartheta^1 \in (\vartheta^* , \vartheta)$
such that $t< T (\vartheta^1, \vartheta^*)$. Next take $\vartheta^2,
\vartheta^3\in \mathds{R}$ such that $\vartheta^1<\vartheta^2 <
\vartheta^3$ and $t< \bar{T} (\vartheta^3, \vartheta^2)$. The latter
is possible since $\bar{T}$ depends only on the difference
$\vartheta_3 - \vartheta_2$, see (\ref{59}). For the fixed $t$, $k_t
\in \mathcal{K}_{\vartheta^1}^\star \hookrightarrow
\mathcal{K}_{\vartheta^3}^\star$, and hence one can write
\begin{eqnarray}
 \label{83}
 u_t & = &
\bar{S}_{\vartheta^3 \vartheta^*}(t) u_0 \\[.2cm]
 & = & (u_0 - k_0)+ k_t + \int_0^t \bar{S}_{\vartheta^3 \vartheta^2} (t-s)
 D_{\vartheta^2\vartheta^1} k_s ds, \nonumber
\end{eqnarray}
where
\begin{equation*}
  D_{\vartheta \vartheta''} = \bar{L}_{\vartheta \vartheta''} -
L^\Delta_{\vartheta \vartheta''}, \qquad    D_\vartheta =
\bar{L}_\vartheta - L^\Delta_\vartheta,
\end{equation*}
and the latter two operators are as in (\ref{60}) and (\ref{24c})
respectively. By Lemma \ref{Id1lm}, for $s\leq t$, $k_s \in
\mathcal{K}_{\vartheta^1}^\star$. By (\ref{15}), (\ref{55}), and
Lemma \ref{J1lm} we have that $D_{\vartheta^2\vartheta^1}:
\mathcal{K}_{\vartheta^1}^\star \to \mathcal{K}_{\vartheta^2}^+$.
Then by Lemma \ref{A1lm} the third summand in the second line in
(\ref{83}) is in $\mathcal{K}_{\vartheta^3}^+$ which completes the
proof since $u_0 - k_0$ is also positive.
\end{proof}
\vskip.1cm \noindent {\it Proof of Theorem \ref{1tm}.} According to
Definition \ref{S1df} and Remark \ref{D1rk} the map $[0,+\infty)\ni
t \mapsto k_t \in \mathcal{K}^\star$ is the solution in question if:
(a) $k_t(\emptyset,\emptyset)=1$; (b) for each $t>0$, there exists
$\vartheta''\in \mathds{R}$ such that $k_t \in
\mathcal{K}_{\vartheta''}$ and  $\frac{d}{dt} k_t =
L^\Delta_\vartheta k_t$ for each $\vartheta> \vartheta''$.

Let $k_0$ and $C>0$ be as in the statement of Theorem \ref{1tm}. Set
$\vartheta^*= \log C$. Then, for $\vartheta = \vartheta^* +
\delta(\vartheta^*)$, see (\ref{N1}) and (\ref{N2}), $k_t$ as given
in (\ref{Sol}) is a unique solution of (\ref{24c}) in
$\mathcal{K}_\vartheta$ on the time interval $[0,T(\vartheta,
\vartheta^*))$. By (\ref{15}) we have
\[
\left(\frac{d}{dt} k_t\right)(\emptyset, \emptyset) = (L^\Delta
k_t)(\emptyset, \emptyset) =0,
\]
which yields that $k_t(\emptyset, \emptyset) = k_0(\emptyset,
\emptyset) =1$. By Lemma \ref{Id1lm} $k_t \in
\mathcal{K}_\vartheta^\star$, and hence $k_t$ is the solution in
question for $t< \tau( \vartheta^*)$.
 According to Lemma \ref{J2lm} $k_t$ lies in
$\mathcal{K}_{\vartheta(t)}$ with $\vartheta (t)$ given in
(\ref{De}). Fix any $\epsilon \in (0,1)$ and then set $s_0=0$, $s_1
=(1-\epsilon) \tau(\vartheta^*)$, and $\vartheta^{*}_1 = \vartheta
(s_1)$. Thereafter, set  $\vartheta^1 = \vartheta^{*}_{1} +
\delta(\vartheta^{*}_{1})$ and
\begin{equation*}
 k_{t+  s_1} = S^1_{\vartheta^1 \vartheta^{*}_{1}}(t)k_{s_1}, \qquad t \in [0, \tau( \vartheta^{*}_{1})).
\end{equation*}
Note that for $t$ such that $t+s_1 < \tau(\vartheta^*)$,
\[
 k_{t+  s_1} = S^1_{\vartheta^1 \vartheta^{*}}(t+s_1)k_{0},
\]
see (\ref{D1}). Thus, by Lemmas \ref{Id1lm} and \ref{J2lm} the map
$[0, s_1 + \tau(\vartheta^{*}_{1})) \ni t \mapsto k_t \in
\mathcal{K}_{\vartheta(t)}$ with
\begin{equation*}
k_t = \left\{ \begin{array}{ll} S^1_{\vartheta^{*}_1
\vartheta^{*}}(t) k_0
\quad &t\leq s_1;\\[.3cm] S^1_{\vartheta^1 \vartheta^{*}_{1}}(t-s_1)
k_{s_1} \quad &t\in [s_1, s_1+ \tau( \vartheta^{*}_{1}))
\end{array} \right.
\end{equation*}
is the solution in question on the indicated time interval. We
continue this procedure by setting $s_n =(1-\epsilon)
\tau(\vartheta^*_{n-1})$, $n\geq 2$, and then
\begin{equation}
  \label{N4}
\vartheta^*_n = \vartheta(s_1+\cdots +s_n ), \qquad \vartheta^n =
\vartheta^*_n + \delta(\vartheta^*_n).
\end{equation}
This yields the solution in question on the time interval $[0, s_1 +
\cdots + s_{n+1}]$ which for $t \in [s_1 +\cdots + s_l, s_1 +\cdots
+ s_{l+1}]$, $l=0, \dots , n$, is given by
\[
k_t =  S^1_{\vartheta^l \vartheta^*_{l}}(t-(s_1 +\cdots +s_l))
k_{s_l}.
\]
Then the global solution in question exists whenever the series
\[
\sum_{n\geq 1}s_n = (1-\epsilon) \sum_{n\geq 1} \tau(\vartheta^*_n)
\]
diverges. Assume that this is not the case. Then by (\ref{De}) and
(\ref{N4}) we get that both (a) and (b) ought to be true, where (a)
$\sup_{n\geq 1} \vartheta^*_n =: \bar{\vartheta}<+ \infty$ and (b)
$\tau(\vartheta^*_n) \to 0$ as $n\to +\infty$. However, by
(\ref{N1}) and (\ref{N2}) it follows that (a) implies
$\tau(\vartheta^*_n) \geq \tau(\bar{\vartheta}) > 0$, which
contradicts (b). \hskip2.5cm $\square$

\section{The Proof of Theorems \ref{Ktm} and \ref{4tm}}
\label{S6}

\subsection{The kinetic equations}
Here we prove Theorem \ref{Ktm}. For a continuous function
\[
[0,+\infty) \ni t \mapsto \varrho_t = (\varrho_{0,t}, \varrho_{1,t})
\in L^\infty (\mathds{R}^d \to \mathds{R}^2),
\]
cf. (\ref{40}), let us consider
\begin{eqnarray}
  \label{41}
& & F_{0,t}(\varrho)(x) = \varrho_{0,0}(x) e^{-\alpha_0 t} \\[.2cm] & & \quad \qquad +
\int_0^t e^{-\alpha_0(t-s)} (a_0 \ast \varrho_{0,s})(x)\exp\left[-
(\phi_0 \ast \varrho_{1,s})(x)\right] ds \nonumber \\[.2cm]
& & \quad \qquad + \int_0^t e^{-\alpha_0(t-s)} \varrho_{0,s}(x)
\bigg{(}a_0 \ast \bigg{[} 1- \exp\left[- (\phi_0 \ast
\varrho_{1,s})\right] \bigg{]} \bigg{)}(x) ds ,\nonumber \\[.2cm]
& & F_{1,t}(\varrho)(x) = \varrho_{1,0}(x) e^{-\alpha_1 t} \nonumber \\[.2cm] & & \quad \qquad +
\int_0^t e^{-\alpha_1(t-s)} (a_1 \ast \varrho_{1,s})(x)\exp\left[-
(\phi_1 \ast \varrho_{0,s})(x)\right] ds \nonumber \\[.2cm]
& & \quad \qquad + \int_0^t e^{-\alpha_1(t-s)} \varrho_{1,s}(x)
\bigg{(}a_1 \ast \bigg{[} 1- \exp\left[- (\phi_1 \ast
\varrho_{0,s})\right] \bigg{]} \bigg{)}(x) ds. \nonumber
\end{eqnarray}
For a given $T>0$, let $\mathcal{C}_T$ stand for the Banach space of
continuous functions
\begin{equation}
  \label{41a}
[0, T] \ni t \mapsto (\varrho_{0,t}, \varrho_{1,t}) \in L^\infty
(\mathds{R}^d \to \mathds{R}^2),
\end{equation}
with norm
\begin{equation}
  \label{42}
\|\varrho\|_T = \max_{i=0,1} \sup_{t\in [0,T]}\left\{
\|\varrho_{i,t}\|_{L^\infty} e^{-\alpha_i t}\right\}.
\end{equation}
Let also $\mathcal{C}_T^+$ denote the set of all positive
$\varrho\in \mathcal{C}_T$, i.e., such that $\varrho_{i,t}(x) \geq
0$ for all $i=0,1$, $t\in [0,T]$, and Lebesgue-almost all $x$. By
means of $F_{i,t}$ introduced in (\ref{41}) we then define the map
\begin{equation*}
\mathcal{C}_T \ni \varrho \mapsto F(\varrho) = (F_0(\varrho),
F_{1}(\varrho)) \in \mathcal{C}_T
\end{equation*}
such that the values of $F_i(\varrho)$ are given in the right-hand
sides of (\ref{41}). By direct inspection one concludes that both
$F_{i,t}(\varrho)$, $i=0,1$, are continuously differentiable in $t$,
and the function as in (\ref{41a}) is a positive solution of
(\ref{29}) on $[0,T]$ if and only if it solves in $\mathcal{C}_T^+$
the following fixed-point equation
\begin{equation}
  \label{44}
\varrho = F(\varrho).
\end{equation}
Let $C>0$ be an arbitrary number and $\varrho_{i,0}$, $i=0,1$, be as
in (\ref{40a}) and (\ref{41}). Set
\begin{equation}
  \label{45}
\varDelta_C = \{ \varrho\in \mathcal{C}_T^+ : (\varrho_{0,t},
\varrho_{1,t})|_{t=0}=(\varrho_{0,0}, \varrho_{1,0}), \ \ {\rm and }
\ \  \|\varrho\|_T \leq C\}.
\end{equation}
By (\ref{41}) one readily gets that $F:\mathcal{C}_T^+\to
\mathcal{C}_T^+$. Let us show that
\begin{equation}
  \label{45a}
\forall C>0 \qquad F:\varDelta_C\to \varDelta_C.
\end{equation}
For $\varrho \in \varDelta_C$, from the first equation in (\ref{41})
one gets
\begin{eqnarray}
  \label{46}
\| F_{0,t} (\varrho) \|_{L^\infty} & \leq &  C e^{-\alpha_0 t} + 2
\alpha_0 e^{-\alpha_0 t} \int_0^t e^{\alpha_0 s} \|\varrho_{0,s}
\|_{L^\infty} d s \\[.2cm] & \leq &  C e^{\alpha_0 t}, \qquad \qquad t \in
[0,T]. \nonumber
\end{eqnarray}
Similarly,  $\| F_{1,t} (\varrho) \|_{L^\infty} \leq C e^{\alpha_1
t}$, which proves (\ref{45a}). To solve (\ref{44}) we apply the
Banach contraction principle. To this end we pick $T>0$ such that
$F$ is a contraction on (\ref{45}). We do this as follows. For
$\varrho, \bar{\varrho} \in \varDelta_C$, like in (\ref{46}) we
obtain
\begin{eqnarray*}
& & \| F_{0,t} (\varrho)  - F_{0,t} (\bar{\varrho}) \|_{L^\infty}
\leq  2 \alpha_0 e^{-\alpha_0 t} \int_0^t e^{\alpha_0 s}
\|\varrho_{0,s} -
\bar{\varrho}_{0,s}\|_{L^\infty} d s \\[.2cm] & & \quad \qquad +
2 \alpha_0 e^{-\alpha_0 t} \int_0^t e^{\alpha_0 s}
\|\bar{\varrho}_{0,s}\|_{L^\infty} \|\varrho_{1,s} -
\bar{\varrho}_{1,s}\|_{L^\infty} d s. \nonumber \\[.2cm] & & \quad \qquad
\leq e^{\alpha_0 t} \|\varrho - \bar{\varrho}\|_T \bigg{(} 1 - e
^{-2 \alpha_0 t}\left[ 1 - \frac{2}{3} C \left( e^{3 \alpha_0 t} -1
\right) \right]\bigg{)} .\nonumber
\end{eqnarray*}
The corresponding estimate for $\| F_{1,t} (\varrho)  - F_{1,t}
(\bar{\varrho}) \|_{L^\infty}$ (with $e^{\alpha_1 t}$) can be
obtained in the same way. Then according to (\ref{42}) $F$ is a
contraction on $\varDelta_C$ whenever $C>0$ and $T$ satisfy
\begin{equation}
  \label{48}
 e^{3 \alpha T} < 1 + \frac{3}{2C}, \quad \qquad \alpha
 :=\max_{i=0,1}\alpha_i.
\end{equation}
This yields the existence of the unique positive solution of
(\ref{29}) on the time interval $[0,T]$, where $T$ is defined in
(\ref{48}) by the initial condition $(\varrho_{0,0},
\varrho_{1,0})$. This solution lies in $\varDelta_C$ and hence
\begin{equation}
  \label{49}
\|\varrho_{i,T}\|_{L^\infty} \leq e^{\alpha T} C, \qquad i=0,1.
\end{equation}
Now we consider the problem (\ref{29}) for $\varrho^{(1)}_{i, t} =
\varrho_{i,T+t}$, $i=0,1$, where $\varrho$ is the solution just
constructed. For this new problem, by (\ref{49}) we have
\[
\|\varrho^{(1)}_{i,0}\|_{L^\infty} \leq C_1 := e^{\alpha T} C,
\qquad i=0,1.
\]
Then we repeat the above construction and obtain the solution
$\varrho^{(1)}$ on the time interval $[0,T_1]$ with $T_1>0$
satisfying, cf. (\ref{48}),
\[
e^{3 \alpha T_1} = 1 + \frac{1}{C} e^{-\alpha T} < 1 +
\frac{3}{2C}e^{-\alpha T} = 1 + \frac{3}{2C_1}.
\]
By further repeating this construction we obtain
$\varrho_{i,t}^{(n)} = \varrho_{i, T+T_1 + \cdots + T_{n-1} + t}$,
$i=0,1$, $t\in [0,T_n]$, where the sequence $\{T_n\}_{n\in
\mathds{N}}$ is defined recursively by the condition
\begin{equation}
  \label{50}
e^{3 \alpha T_n} = 1 +  \frac{1}{C}\exp\left[ -\alpha \left(T + T_1
+\cdots + T_{n-1}\right) \right], \quad n\in \mathds{N}.
\end{equation}
Thus, the global solution in question exists if the series $\sum_{n}
T_n$ is divergent. Assume that this is not the case. Then the
right-hand side of (\ref{50}) is bounded from below by some $b
>1$, uniformly in $n$. This yields that $T_n \geq \log b/ 3
\alpha >0$, holding for all $n\in \mathds{N}$, which contradicts the
summability of $\{T_n\}_{n\in \mathds{N}}$ and thus completes the
proof of Theorem \ref{Ktm}.

\subsection{The scaling limit}

For each $k$ and $\lambda$-almost all $\eta \in \Gamma^2_0$, we have
that the following holds, cf. (\ref{13}) and (\ref{31a}),
\begin{eqnarray*}
& & (Q^0_{y,\varepsilon} k) (\eta_0,\eta_1) \to (Q^0_{y,0} k)
(\eta_0,\eta_1)\\[.2cm] & & \qquad \qquad := \int_{\Gamma_0} k(\eta_0, \eta_1\cup \xi)
e(-\phi_0 (y - \cdot);\xi)\lambda (d\xi), \quad \varepsilon \to 0,\nonumber \\[.2cm]
& & (Q^1_{y,\varepsilon} k) (\eta_0,\eta_1) \to (Q^1_{y,0} k)
(\eta_0,\eta_1) \nonumber \\[.2cm] & & \qquad \qquad := \int_{\Gamma_0} k(\eta_0\cup \xi, \eta_1)
e(-\phi_1 (y - \cdot);\xi)\lambda (d\xi), \quad \varepsilon \to 0.
\nonumber
\end{eqnarray*}
Thus, for each $k$ and $\lambda$-almost all $\eta \in \Gamma^2_0$,
\[
(L^{\varepsilon, \Delta} k)(\eta) \to (Vk)(\eta), \qquad {\rm as } \
\ \varepsilon \to 0,
\]
where, cf (\ref{15}),(\ref{31a})
\begin{eqnarray}
  \label{35}
 (Vk)(\eta_0 , \eta_1)  & = & \sum_{y\in \eta_0}
\int_{\mathds{R}^d} a_0 (x-y) (Q_{y,0}^0 k) (\eta_0\setminus y \cup
x,
\eta_1) d x \nonumber \\[.2cm]
& - & \sum_{x\in \eta_0} \int_{\mathds{R}^d} a_0 (x-y) (Q_{y,0}^0 k)
(\eta_0,
\eta_1) d y  \\[.2cm]
& + & \sum_{y\in \eta_1} \int_{\mathds{R}^d} a_1 (x-y) (Q_{y,0}^1 k)
(\eta_0,
\eta_1\setminus y \cup x) d x \nonumber \\[.2cm] & - & \sum_{x\in \eta_1} \int_{\mathds{R}^d} a_1 (x-y)
(Q_{y,0}^1 k) (\eta_0, \eta_1) d y. \nonumber
\end{eqnarray}
Like above, for each $\vartheta''\in \mathds{R}$ and $k\in
\mathcal{K}_{\vartheta''}$, both $Q^i_{y,0}k$ satisfy the estimates
as in (\ref{21}) and (\ref{22}). Then for $\vartheta,\vartheta''\in
\mathds{R}$ such that $\vartheta'' < \vartheta$, $\|V
k\|_{\vartheta}$ is bounded by the right-hand side of (\ref{24}).
This allows one to define the operators $V_\vartheta$ and
$V_{\vartheta\vartheta''}$ analogous to $L^\Delta_\vartheta$ and
$L^\Delta_{\vartheta\vartheta''}$, respectively. For $\varrho_t$
being the solution as in Theorem \ref{Ktm}, $k_{\pi_{\varrho_t}}$
satisfies
\begin{equation}
  \label{J}
\frac{d}{dt} k_{\pi_{\varrho_t}} = V_{\vartheta \vartheta''}
k_{\pi_{\varrho_t}}, \qquad t>0,
\end{equation}
where $\vartheta''\in \mathds{R}$ is such that
$k_{\pi_{\varrho_t}}\in \mathcal{K}_{\vartheta''}$, see (\ref{40a}),
and $\vartheta > \vartheta''$ is arbitrary.  This can be checked by
direct calculations based on (\ref{35}) and (\ref{29}). Moreover, if
we set $C= \|\varrho_0\|_\infty$, see (\ref{l-norm}), then
$k_{\pi_{\varrho_t}}$ satisfies (\ref{24d}) with this $C$, which
follows from (\ref{40a}).
Thus, by Corollary \ref{D1co} we conclude
that $k_{\pi_{\varrho_t}}\in \mathcal{K}_{\vartheta(t)}$ for all
$t>0$. \vskip.1cm \noindent
 \textit{The proof of Theorem \ref{4tm}.}
Let $\vartheta_*$ be as assumed. As mentioned above, we then have
that $k_{\pi_{\varrho_t}} \in \mathcal{K}_{\vartheta_T}$ for all
$t\in [0,T]$ with $\vartheta_T := \vartheta_* + \alpha T$ and $T$
such that
\begin{equation}
  \label{J0}
 T < \tau(\vartheta_* + \alpha T).
\end{equation}
The latter is possible since the function $\vartheta\mapsto
\tau(\vartheta)$ is continuous and $\tau(\vartheta_* ) >0$, see
(\ref{N1}). Since the inequality in (\ref{J0}) is strict, we can
also pick $\vartheta_1 > \vartheta_T$ such that $T <
\tau(\vartheta_1)$. Thereafter, we set $\vartheta = \vartheta_1 +
\delta (\vartheta_1)$, cf. Remark \ref{JJrk}. For $q_{0,
\varepsilon}$ with the assumed property, let $q_{t,\varepsilon}$ be
the solution of (\ref{32a}) in $\mathcal{K}_\vartheta$. In view of
(\ref{J}), we then have
\begin{gather}
  \label{J0a}
q_{t,\varepsilon} - k_{\pi_{\varrho_t}} = \int_0^t
S^\varepsilon_{\vartheta \vartheta_1}(t-s) \left(L^{\varepsilon,
\Delta}_{\vartheta_1 \vartheta_T} - V_{\vartheta_1 \vartheta_T}
\right)
k_{\pi_{\varrho_s}} d s \\[.2cm] + S^\varepsilon_{\vartheta
\vartheta_*}(t)\left[ q_{0, \varepsilon} -
k_{\pi_{\varrho_0}}\right], \qquad t \in [0,T]. \nonumber
\end{gather}
Since $\vartheta \mapsto \tau(\vartheta)$ is decreasing, by
(\ref{J0}) we have that $T < \tau(\vartheta_*)$. By (\ref{51}) we
then get
\[
\forall t \in [0,T] \qquad  \|S^\varepsilon_{\vartheta
\vartheta_*}(t) \| \leq \frac{T(\vartheta,\vartheta_*)}{
T(\vartheta,\vartheta_*) - T},
\]
which yields that the second term in (\ref{J0a}) tends to zero
uniformly on $[0,T]$. Also by (\ref{51}) we have
\begin{eqnarray}
  \label{J0b}
& & \bigg{\|} \int_0^t S^\varepsilon_{\vartheta \vartheta_1}(t-s)
\left(L^{\varepsilon, \Delta}_{\vartheta_1 \vartheta_T} -
V_{\vartheta_1 \vartheta_T} \right) k_{\pi_{\varrho_s}} d s
\bigg{\|}_\vartheta  \\[.2cm] & & \qquad \leq \|k_{\pi_{\varrho_T}}\|_{\vartheta_T} \tau(\vartheta_1) \log
\frac{T(\vartheta,\vartheta_1)}{T(\vartheta,\vartheta_1) -T} \|
L^{\varepsilon, \Delta}_{\vartheta_1 \vartheta_T} - V_{\vartheta_1
\vartheta_T}\|.\nonumber
\end{eqnarray}
To estimate the latter term we set
\begin{equation}
  \label{86}
W^i_{y,\varepsilon}k = Q^i_{y,0}k - Q^i_{y,\varepsilon}k, \qquad
i=0,1, \ \ y\in \mathds{R}^d.
\end{equation}
By means of the inequality, cf. the proof of Theorem 4.6 in
\cite{BKKK},
\[
\left\vert b_1 \cdots b_n - a_1 \cdots a_n \right\vert \leq
\sum_{i=1}^n |b_i - a_i| \prod_{j\neq i} \max\{|a_j|; |b_j|\},
\]
and
\begin{equation*}
0\leq \psi(t) := (t -1 +e^{-t})/t^2 \leq 1/2, \qquad t\geq 0,
\end{equation*}
 we obtain, cf. (\ref{21}),
\begin{eqnarray}
  \label{88}
& & \left\vert W^0_{y,\varepsilon}k (\eta_0, \eta_1)\right\vert \leq
\varepsilon \|k\|_{\vartheta''} \exp\left( \vartheta'' |\eta_0| +
\vartheta''
|\eta_1| \right) \\[.2cm]
&  \times &  \int_{\Gamma_0}\bigg{(} e^{\vartheta'' |\xi|}
\sum_{x\in \xi} \left[ \phi_0(y-x) \right]^2 \psi \left( \varepsilon
\phi_0(y-x)\right) \prod_{z\in \xi\setminus x} \phi_0(y-z)
\bigg{)}\lambda(d \xi)
 \nonumber  \qquad
 \\[.2cm]
&  \leq &( \varepsilon /2)\bar{\phi_0} \|k\|_{\vartheta''}
\exp\left(\vartheta'' |\eta_0| + \vartheta'' |\eta_1| \right)
\int_{\Gamma_0}\bigg{(} |\xi| e^{\vartheta'' |\xi|} \prod_{z\in \xi}
\phi_0(y-z) \bigg{)}\lambda(d \xi) \nonumber \\[.2cm]
& = & (\varepsilon/2) \bar{\phi_0}  \langle \phi_0 \rangle
\exp\left( \langle \phi_0 \rangle e^{\vartheta''} \right)
\|k\|_{\vartheta''} \exp\left( \vartheta'' |\eta_0| + \vartheta''
|\eta_1| +\vartheta'' \right). \nonumber
\end{eqnarray}
Likewise,
\begin{eqnarray}
  \label{88a}
& & \left\vert W^1_{y,\varepsilon}k (\eta_0, \eta_1)\right\vert \leq
\\[.2cm]& & \qquad (\varepsilon/2) \bar{\phi_1}  \langle \phi_1 \rangle \exp\left(
\langle \phi_1 \rangle e^{\vartheta''} \right) \|k\|_{\vartheta''}
\exp\left( \vartheta'' |\eta_0| + \vartheta'' |\eta_1| +\vartheta''
\right). \nonumber
\end{eqnarray}
Next, by (\ref{15}), (\ref{35}), and (\ref{86}) we have
\begin{eqnarray}
  \label{350}
 & & (L^{\varepsilon,\Delta} - V)k(\eta_0 , \eta_1)   =  \sum_{y\in \eta_0}
\int_{\mathds{R}^d} a_0 (x-y) (U^0_{y,\varepsilon} k)
(\eta_0\setminus y \cup x,
\eta_1) d x \qquad \qquad  \nonumber \\[.2cm]
& & \qquad \quad -  \sum_{x\in \eta_0} \int_{\mathds{R}^d} a_0 (x-y)
(U^0_{y,\varepsilon} k) (\eta_0,
\eta_1) d y  \\[.2cm]
& & \qquad \quad  +  \sum_{y\in \eta_1} \int_{\mathds{R}^d} a_1
(x-y) (U^1_{y,\varepsilon} k) (\eta_0,
\eta_1\setminus y \cup x) d x \nonumber \\[.2cm]& & \qquad \quad  -  \sum_{x\in \eta_1} \int_{\mathds{R}^d} a_1 (x-y)
(U^1_{y,\varepsilon} k) (\eta_0, \eta_1) d y. \nonumber
\end{eqnarray}
Here we use the following notations
\begin{gather*}
(U^0_{y,\varepsilon} k) (\eta_0,\eta_1) =
e(\tau^0_{y,\varepsilon};\eta_1) (Q^0_{y,\varepsilon} k)
(\eta_0,\eta_1) - (Q^0_{y,0} k) (\eta_0,\eta_1), \\[.2cm]
(U^1_{y,\varepsilon} k) (\eta_0,\eta_1) =
e(\tau^1_{y,\varepsilon};\eta_0) (Q^1_{y,\varepsilon} k)
(\eta_0,\eta_1) - (Q^1_{y,0} k) (\eta_0,\eta_1). \nonumber
\end{gather*}
Then, cf. (\ref{86}),
\begin{eqnarray*}
\left\vert (U^0_{y,\varepsilon} k) (\eta_0,\eta_1)\right\vert & \leq
& \left\vert (W^0_{y,\varepsilon} k) (\eta_0,\eta_1)\right\vert +
\varepsilon \bar{\phi}_0 |\eta_1|\left\vert(Q^0_{y,0} k)
(\eta_0,\eta_1)\right\vert, \qquad \\[.2cm]
\left\vert (U^1_{y,\varepsilon} k) (\eta_0,\eta_1)\right\vert & \leq
& \left\vert (W^1_{y,\varepsilon} k) (\eta_0,\eta_1)\right\vert +
\varepsilon \bar{\phi}_1 |\eta_0|\left\vert(Q^1_{y,0} k)
(\eta_0,\eta_1)\right\vert.\nonumber
\end{eqnarray*}
Now by (\ref{21}), (\ref{22}), (\ref{88}), (\ref{88a}), and
(\ref{350}) we get
\begin{eqnarray*}
  \label{353}
 & & \left\vert(L^{\varepsilon,\Delta} - V)k(\eta_0 ,
 \eta_1)\right\vert \leq
\varepsilon \alpha \|k\|_{\vartheta''} \exp\left[ \vartheta''
(|\eta_0|+|\eta_1|)\right]\exp\left( ce^{\vartheta''}\right) \qquad
\\[.2cm]
& & \qquad \times \bigg{(}
2|\eta_0||\eta_1|\left(\bar{\phi}_0+\bar{\phi}_1\right) +
e^{\vartheta''} \left(\bar{\phi}_0 \langle \phi_0 \rangle |\eta_0| +
\bar{\phi}_1 \langle \phi_1 \rangle |\eta_1| \right) \bigg{)}
\end{eqnarray*}
Like in obtaining (\ref{24}) we then get from the latter
\begin{equation*}
\|L^{\varepsilon,\Delta}_{\vartheta_1 \vartheta_T} - V_{\vartheta_1
\vartheta_T}\| \leq \varepsilon \varPhi(\vartheta_1, \vartheta_T),
\end{equation*}
where $\varPhi(\vartheta_1, \vartheta_T)>0$ depends  on the choice
of $\vartheta_1$, $\vartheta_T$ and on the model parameters only,
and may be calculated explicitly. Then the use of the latter
estimate in (\ref{J0b}) completes the proof.

\hfill%
$\square $

\section*{Acknowledgment}
The present research was supported by the European Commission under
the project STREVCOMS PIRSES-2013-612669.

\end{document}